\documentclass{article}

\usepackage[utf8]{inputenc}
\usepackage[T1]{fontenc}
\usepackage{amsmath,amssymb,amsthm,bbm,float,graphicx,geometry,lmodern,mathtools,parskip,setspace,subcaption}
\usepackage[colorlinks, linkcolor = blue!80!black, citecolor = blue!80!black,breaklinks, pdfauthor={Benedikt Jahnel, Lukas Luechtrath, Anh Duc Vu}]{hyperref}
\usepackage{mathrsfs}
\usepackage{enumerate}
\usepackage{nicefrac}
\usepackage{todonotes}
\usepackage{comment}
\usepackage{titling}
\usepackage{fullpage}
\usepackage{comment}
\usepackage{orcidlink}
\usepackage{framed}
\usepackage{xcolor}
\usepackage[most]{tcolorbox}
\usepackage{xparse} % for optional arguments
\usepackage{dsfont}
\usepackage[shortlabels]{enumitem}

\usepackage{tikz}
\usetikzlibrary{backgrounds}
\usetikzlibrary{patterns}
\usetikzlibrary{positioning, shapes.geometric}
\usetikzlibrary{external, arrows, decorations.pathmorphing} 

\usepackage{caption}
\usepackage{graphicx}
\graphicspath{{Images/}}%put any images in this file
\allowdisplaybreaks%breaks long equations onto multiple pages is needed
\usepackage[raggedright]{titlesec} % Takes care of too long section titles 
\usepackage{booktabs} %% Nice tables
\newcommand{\N}{\mathbb{N}}

\DeclareMathOperator{\interior}{int}

\renewcommand*{\P}{\mathbb{P}}
\newcommand*{\E}{\mathbb{E}}
\newcommand*{\R}{\mathbb{R}}

\newcommand*{\Z}{\mathbb{Z}}

\renewcommand{\N}{\mathbb{N}}
\newcommand*{\X}{\mathbb{X}}
\newcommand*{\Y}{\mathbb{Y}}

\newcommand*{\e}{\mathrm{e}}
\renewcommand*{\Xi}{\varXi}
\renewcommand*{\epsilon}{\varepsilon}
\renewcommand*{\theta}{\vartheta}
\renewcommand*{\Theta}{\varTheta}

\renewcommand*{\Delta}{\varDelta}

\newcommand{\Zz}{\mathbb{Z}_{\geq0}}
\newcommand{\LRC}{\mathrm{LRC}}
\newcommand{\RLC}{\mathrm{RLC}}
\newcommand{\DTC}{\mathrm{BTC}}
\newcommand{\hor}{\mathrm{hor}}
\newcommand{\ver}{\mathrm{ver}}
\newcommand{\Ltwo}{\mathbb{L}^2}
\newcommand{\PPP}{\mathrm{PPP}}

\newcommand{\eps}{\varepsilon}

\renewcommand{\interior}[1]{%
  {\kern0pt#1}^{\mathrm{o}}%
}

\newtheorem{theorem}{theorem}[section]
\newtheorem{cor}[theorem]{Corollary}
\newtheorem{lemma}[theorem]{Lemma}
\newtheorem{proposition}[theorem]{Proposition}
\newtheorem{thm}[theorem]{Theorem}
\theoremstyle{definition}
\newtheorem{definition}[theorem]{\bf Definition}

\theoremstyle{remark}

\usepackage{thmtools}
\usepackage{etoolbox}

\declaretheoremstyle[
  shaded={rulecolor=black, rulewidth=1pt, bgcolor=gray!10},
  name=Parameters,
  headfont=\bfseries,
  sibling = theorem,
]{stypar}

\declaretheorem[style=stypar]{parameters}

\usepackage{nameref}

\makeatletter
\let\orgdescriptionlabel\descriptionlabel
\renewcommand*{\descriptionlabel}[1]{%
  \let\orglabel\label
  \let\label\@gobble
  \phantomsection
  \edef\@currentlabel{#1}%
  \let\label\orglabel
  \orgdescriptionlabel{#1}%
}

\pretitle{\centering\LARGE\scshape}
 \posttitle{\vskip 0.75cm}

 \predate{\vskip 0.75 cm \centering\large}
 \postdate{\par}

\usepackage[style = numeric, sorting=nyt, url = false, abbreviate=false, maxbibnames=9, sortcites=true, doi = true, backend = biber, giveninits = true, isbn=false]{biblatex}
\renewbibmacro{in:}{\ifentrytype{article}{}{\printtext{\bibstring{in}\intitlepunct}}}
\bibliography{references_2}

\title{Oriented bond-site percolation in random environment and contact processes with periodic recovery}
\thanksmarkseries{arabic}

\author{
Benedikt Jahnel 
	\orcidlink{0000-0002-4212-0065}
	\thanksgap{0.4ex}
	\thanks{Technische Universit\"at Braunschweig, Universit\"atsplatz 2, 38106 Braunschweig, Germany} 
	\thanksgap{0.4ex}
	\thanks{Weierstrass Institute for Applied Analysis and Stochastics, Mohrenstraße 39, 10117 Berlin, Germany} 
	\\ benedikt.jahnel@tu-braunschweig.de 
\and
Lukas L\"{u}chtrath 
	\orcidlink{0000-0003-4969-806X} 
	\thanksmark{2}
	\\ lukas.luechtrath@wias-berlin.de
\and
\\
Anh Duc Vu 
	\orcidlink{0009-0005-6913-4992} 
	\thanksmark{2}
	\\ anhduc.vu@wias-berlin.de
}

\date{\today}

\usepackage{cleveref}

\begin{document}

\maketitle

\begin{spacing}{0.9}
\begin{abstract} 
\noindent 
We investigate oriented bond-site percolation on the planar lattice in which entire columns are stretched. Generalising recent results by Hilário et al., we establish non-trivial percolation under a $(1+\eps)$-th moment condition on the stretches and use this to prove survival of contact processes with periodic recoveries as well as in random environments.

\medskip
\noindent\footnotesize{{\textbf{AMS-MSC 2020}: Primary: 60K05; Secondary: 60K35, 82B43}

}
\end{abstract}
\end{spacing}

%%%%%%%%%%%%%%%%%%%%%%%%%%%%%%%%%%%%%%%%%%%%%%%%%%%%%%%%%%%%%%%%%%%%%%%%%%%%%%%%%%%%%%%%%%%%%%%
%%%%%%%%%%%%%%%%%%%%%%%%%%%%%%%%%%%%%%%%%%%%%%%%%%%%%%%%%%%%%%%%%%%%%%%%%%%%%%%%%%%%%%%%%%%%%%%
%%%%       Introduction
%%%%%%%%%%%%%%%%%%%%%%%%%%%%%%%%%%%%%%%%%%%%%%%%%%%%%%%%%%%%%%%%%%%%%%%%%%%%%%%%%%%%%%%%%%%%%%%
%%%%%%%%%%%%%%%%%%%%%%%%%%%%%%%%%%%%%%%%%%%%%%%%%%%%%%%%%%%%%%%%%%%%%%%%%%%%%%%%%%%%%%%%%%%%%%%

\section{Setting, results and discussion}\label{sec:intro}
Our motivation stems from the desire to understand survival in generalised contact processes, in particular in random spatial environments and in situations where the Poisson point processes in the graphical construction are replaced by much more rigid point processes such as randomly shifted lattices, see Section~\ref{subsec:Application-CP}. In order to guarantee survival, we rely on the analysis of certain stretched oriented bond-site percolation models, see Section~\ref{subsec:OPRE}. For these, we establish the existence of percolation phases, which constitutes our main result, Theorem~\ref{thm:OPRE-percolates}. 

\subsection{Oriented percolation}\label{subsec:OPRE}
\begin{figure}
    \centering  \includegraphics[width=0.95\linewidth]{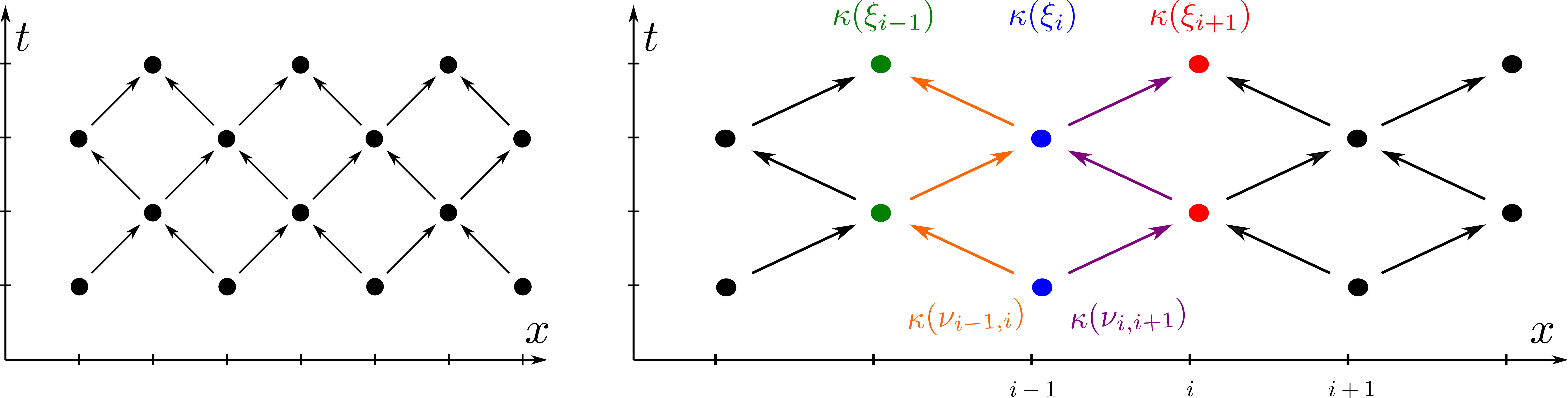}
    \caption{Illustrations of the $\Ltwo$-lattice (left) and the OPRE on $\Ltwo$ (right). Coloured vertices and edges of the same colour have the same probability to be open; for example, blue vertices are open with probability $\kappa(\xi_x)$.}
    \label{fig:OPRE-and-L2-lattice}
\end{figure}

Consider the oriented $\Ltwo$-lattice, i.e., the graph $\Ltwo=(V,E)$ defined by
\begin{equation}\label{eq:L2-lattice}
    V=\{(t,x)\in\Zz^2 \colon t+x \text{ even}\} 
    \qquad\text{and}\qquad 
    E=\{(t,x)\to(t+1,x+z) \colon (t,x)\in V,\,z\in\{-1,1\}\}.
\end{equation}
We interpret the second component as {\em space} and the first component as {\em time}, noting that time always increases along (directed) paths. Building on this, we generalise classical i.i.d.\ Bernoulli bond-site percolation by introducing environments that weaken both bonds and sites via stretches. These stretches are associated with the spatial component (i.e., the $x$-coordinate), and we refer to the resulting oriented bond-site percolation model as OPRE, short for \emph{oriented percolation in a random environment}; see Figure~\ref{fig:OPRE-and-L2-lattice} for an illustration.

\begin{definition}[Oriented percolation in (spatial columnar) random environment (OPRE)]
\label{def:OPRE-model}
    Fix some spatial stretches $(\xi_x)_{x\in\Zz}$, $(\nu_{x,x+1})_{x\in\Zz}\subset[0,\infty)$, called the \emph{environment}, and a {\em connection function} $\kappa\colon[0,\infty)\to[0,1]$. Consider the following independent bond-site percolation model on $\Ltwo$ in the environment $(\xi_x)_{x\in\Zz}$, $(\nu_{x,x+1})_{x\in\Zz}$, where 
    \begin{itemize}
    	\item 
    		each vertex $(t,x)\in V$ is open independently with probability \(\kappa(\xi_x)\),
		\item 
			each edge of the form \((t,x)\to(t+1,x+1)\) is open independently with probability \(\kappa(\nu_{x,x+1})\), and
		\item	
			each edge of the form \((t,x)\to(t+1,x-1)\) is open independently with probability \(\kappa(\nu_{x-1,x})\).
    \end{itemize}
    If the graph contains an infinite path consisting of open vertices and open edges only, we say that the model \emph{percolates}. 
    If the environment is given by a realisation of families of random variables, we refer to the associated percolation model as \emph{OPRE}.
\end{definition}

We now present our main theorem, which is an extension of
principal results in~\cite{hilario2023phasetransition}, and establishes percolation in OPRE.

\begin{thm}[Percolation]\label{thm:OPRE-percolates}
    Assume that the environment $(\xi_x)_{x\in\Zz},(\nu_{x,x+1})_{x\in\Zz}$ in Definition~\ref{def:OPRE-model} is given by mutually independent families of i.i.d.\ random variables that satisfy the moment condition 
    \begin{equation}\label{eq:OPRE-stretches-moment-condition}
        \E[\xi_0^{1+\eps}], \, \E[\nu_{0,1}^{1+\eps}]<\infty\qquad\text{for some $\eps>0$}.
    \end{equation}
 Furthermore, let $(\kappa_\lambda)_{\lambda\geq 0}$ be a family of connection functions satisfying
    \begin{enumerate}[(i)]
        \item $\kappa_\lambda\colon [0,\infty)\to[0,1]$ is monotonically decreasing for every $\lambda\geq0$, and
        \item there exists some $\sigma>0$ such that, for every $s\in\N$,        \begin{equation}\label{eq:OPRE-requirements-connection-kernel}
          \exp(-\sigma s) \leq \kappa_\lambda(s)  \quad \text{ for every }\lambda\geq0, 
          \qquad\text{ as well as } \qquad
          \kappa_\lambda(s)\xrightarrow{\lambda\to\infty} 1.
        \end{equation} 
    \end{enumerate}    
    Then, there exists $\lambda_c  <\infty$ such that, for every $\lambda>\lambda_c$, the OPRE almost-surely percolates for almost-all realisations of the environment $(\xi_x)_{x\in\Zz},(\nu_{x,x+1})_{x\in\Zz}$.
\end{thm}

Due to ergodicity of the environment, percolation is equivalent to the existence of an infinite open path starting at the origin $(0,0)$ with positive probability. Note further that we may also use two different connection functions for vertices and edges since Theorem~\ref{thm:OPRE-percolates} remains valid when replacing them with their minimum. Although the monotonicity assumption on $\kappa_\lambda$ can be relaxed, we retain it to avoid technical complications. 
A prototypical example is given by $\kappa_{\lambda}(s)=(1-\exp(-\lambda))^s$ and stretches given by exponential random variables. In this setting, the theorem ensures percolation for all sufficiently large $\lambda$. 

\subsection{Related work and outlook}
OPRE, as we present it, can be seen as a generalising framework that includes several models that have been studied previously in the literature. First of all, by choosing $\kappa_\lambda\equiv 1-\exp(-\lambda)$, we obtain classical {\em oriented Bernoulli bond-site percolation} on $\Ltwo$ whereas for $\kappa_\lambda(s)=(1-\exp(-\lambda))^s$, OPRE represents a bond-site version of the {\em randomly stretched lattice} (RSL), originally introduced in~\cite{jonasson2000rsl} on $\Z^d$. 

More precisely, the original RSL is an undirected bond-percolation model in which edges $(x,x+e_i)$ are open independently with probability $p^{\nu_{xe_i}(i)}$, where $e_i$ is the $i$-th standard unit vector, $p\in (0,1)$, and $(\nu_\ell(i))_{\ell\in\Z,\, 1\le i\le d}\subset [0,\infty)$ is a given random environment. This random environment is usually built using i.i.d.\ geometric random variables. The absence of percolation in the RSL for small $p$ follows from a direct coupling to standard Bernoulli bond percolation and the existence of percolation was shown in~\cite{jonasson2000rsl} for $d\geq3$. This has been improved in~\cite{hoffman2005phasetransition}, where percolation was established also for $d=2$, for sufficiently large $p$ and small geometric stretches $(\nu_\ell(i))_{\ell,i}$ and has further been improved to $p>p_c=1/2$ in~\cite{deLima2023dependent} and~\cite{hilario2023newproof}. The latter work also proves that light-tailed stretches are both necessary and sufficient for percolation and we observe  analogous behaviour for a version of OPRE where temporal stretches (in the $t$-coordinate) instead of spatial stretches (in the $x$-coordinate) are used; see Proposition~\ref{prop:OPRE-temporal-extinction} below.
Still in case of the RSL, the situation changes when the planar lattice is only stretched in {\em one} direction. In this case, existence of $(1+\eps)$-th moments are already sufficient for percolation, see~\cite{hilario2023phasetransition}. In fact, this moment condition has been relaxed to $\E[\nu \exp(c(\log\nu)^{1/2})]<\infty$ for some $c>64$ in~\cite{guedes2024phase}. 
%Although our work is based on~\cite{hilario2023phasetransition}, we are confident that our proof modifications are robust enough to also be applicable to the situation considered in~\cite{guedes2024phase}. 
Let us also mention here related research for percolation close to the critical value~\cite{duminilcopin2018brochette}, on a diluted model in $\Z^3$~\cite{hilario2022strictinequality}, for $\ell$-dependent stretches on trees~\cite{moreira2024treesdependent}, and regarding electrical conductivity on a further stretched lattice~\cite{jahnel2024conductance}. 

Returning to OPRE, we note that the main challenges in the analysis of the model are already present in the different versions of RSL just discussed, and consist in controlling the non-decorrelating dependencies induced by the random environment. To tackle these challenges, we use properly adjusted {\em multiscale renormalisation schemes}, inspired by the works mentioned above, in which the environment is treated on different scales and is subdivided into good and bad blocks.

Let us recall that our main motivation to analyse OPRE is to use it as a tool to understand non-standard infection processes and similar systems such as moving populations. A classical model for this is the {\em contact process}, which captures both spatial spread and recovery of individuals. In this manuscript, we treat, for example, inhomogeneities caused by hostile environments, whether spatial or temporal, and demonstrate that the question of survival can be addressed via the existence of infinite open paths in OPRE. To this end, recall the classical Harris contact process, a Markovian model for the spread of infections. Infected vertices transmit the infection along (directed) edges at some fixed rate $\lambda>0$ while recovering independently at a rate $\mu>0$. This dynamic may perhaps be best described through its {\em graphical representation}: Infections are transmitted along an edge at times given by independent Poisson point processes with parameter $\lambda\ge 0$, while recoveries occur at times given by independent Poisson point processes with parameter $\mu \ge 0$. The survival-extinction phase transition of a contact process in a random environment, where the recovery rate $\mu$ depends on the location, was established in~\cite{bramson1991contact}. Given a random environment, the resulting model remains Markovian and can be coupled with an OPRE. This coupling enables us to establish a stronger version of the results in~\cite{bramson1991contact}, as presented in Proposition~\ref{prop:CPRE-BDS}. Let us mention in this context that the phase transition of a related contact-process model containing both spatial and temporal stretches was established in~\cite{jahnel2023contact}, albeit at the cost of requiring long-range edges.

Moreover, OPRE can also be exploited for the analysis of non-standard contact processes in another way. Returning to the graphical representation discussed above, a recent line of research has moved beyond the Markovian framework by replacing the involved Poisson processes by more general point processes, introduced as the {\em generalised contact process} in~\cite{hilario2022contactprocess}.
%, where general point processes are employed in the graphical representation. 
The authors establish extinction in two models, the {\em renewal contact process} and the contact process with dynamic edges, cf.\ also~\cite{linker2020contact}.  In particular, they show that if recovery times follow sufficiently heavy-tailed renewal processes, then the infection survives with positive probability, regardless of the infection rate $\lambda$, see also~\cite{fontes2019contactI}).
 Remarkably, this effect can be so strong that survival is possible even on finite graphs~\cite{fontes2021renewalfinitegraph}. In a follow-up work~\cite{fontes2020contactII}, a subcritical phase is established for sufficiently light-tailed recovery processes, either dimension one or in higher dimensions depending on additional model specifications. The results in~\cite{fontes2019contactI} were further strengthened in~\cite{fontes2023renewal}, where complete convergence for the surviving infection is shown. A more detailed survey can be found in~\cite{fontes2023overviewRCP}. Nonetheless, the general question of survival under light-tailed recovery remains open. A partial answer is provided in~\cite{santos2023survival}, where survival is proven for renewal processes with continuous interarrival distributions having bounded support. 
As an application of our main Theorem~\ref{thm:OPRE-percolates} on OPRE, we address this question for almost {\em deterministic} interarrival times. Observe that, if the recovery times occur exactly at integer times (i.e., recovery times are given by $\mathbb{Z}$), then all vertices recover simultaneously at time 1, causing the infection to die out immediately. However, this argument  this argument fails once a random delay is introduced. In Section~\ref{sec:contact-processes}, we establish the phase transition for survival of the infection under some minimal continuous and discrete random delays in the recovery times.

Finally, we briefly outline potential directions for future research. The existence of the $(1+\eps)$-th moment in Theorem \ref{thm:OPRE-percolates} is used to obtain a decoupling inequality for the random environment \cite{hilario2023phasetransition,guedes2024phase} but it remains an open question whether finite first moments already suffice. It would also be interesting to see the effects of adding a combination of temporal and spatial stretches to the $\Ltwo$-lattice. In the broader context of contact processes and other interacting particle systems with graphical representations, one could, for example, employ random closed sets, as exemplified by first contact percolation in~\cite{jahnel2025fcp}. What general results can be obtained in this more abstract setting? Similarly, we restrict ourselves to periodic recoveries in this paper but the complementary setting of periodic infections remains open. While survival in this setting is typically straightforward to establish, identifying a non-trivial extinction phase appears to be substantially more challenging, especially in higher dimensions.

\subsection{Oriented percolation under temporal stretches}\label{subsec:temporal-stretches}
This section is devoted to rigorously discuss temporal stretches as mentioned above.
%We have mentioned temporal stretches above. Now, we want to discuss this setting rigorously.
Let $(\nu_{t,t+1})_{t\in\Zz}\subset[0,\infty)$ be a family of temporal stretches and $p\in(0,1)$. Consider independent bond percolation on $\Ltwo$, where an edge $(t,x)\to(t+1,x+z)$ with $z\in\{-1,1\}$ is independently open with probability $p^{\nu_{t,t+1}}$. In case of independent geometric distributed temporal stretches, it is known that the oriented percolation model almost-surely contains an infinite open path for all sufficiently large $p$, see~\cite[Thm.~1.1]{kesten2022orientedpercolation} and~\cite[Thm.~8.2]{hilario2023newproof}. The following fact emphasises that temporal stretches behave differently from spatial ones.

\begin{proposition}[Absence of phase transition under large temporal stretches]\label{prop:OPRE-temporal-extinction}
    Assume that $(\nu_{t,t+1})_{t\in\Zz}$ is a family of i.i.d.\ random variables with heavier than exponential tails, i.e.,
    \begin{equation}\label{eq:heavy-tail-stretches}
        \limsup_{s\to\infty} c^s\P(\nu_{0,1} > s) =\infty \qquad\text{ for all } c>1.
    \end{equation}
    Then, for every $p\in(0,1)$, the oriented percolation model with temporal stretches almost-surely does not contain an infinite open path.
\end{proposition}
The proof is given in Section~\ref{subsec:proof-temporal-environment} and uses the fact that the number of reached vertices grows at most polynomially over time in an unstretched model. Hence, the result applies mutatis mutandis in higher dimensions and similar settings, e.g., contact processes where a global recovery rate randomly changes over time with heavier-than-exponential tails.

\subsection{Contact processes with periodic recovery}\label{subsec:Application-CP}
We have already discussed above that survival of contact processes is closely related to oriented percolation. In this section, we first formally introduce the {\em generalised contact process} where infection and recovery times are given by {\em random closed sets} (with respect to the Fell-topology) instead of Poisson processes, extending the notion of the generalised contact process in~\cite{hilario2022contactprocess}.

\begin{figure}
    \centering
    \includegraphics[width=0.8\linewidth]{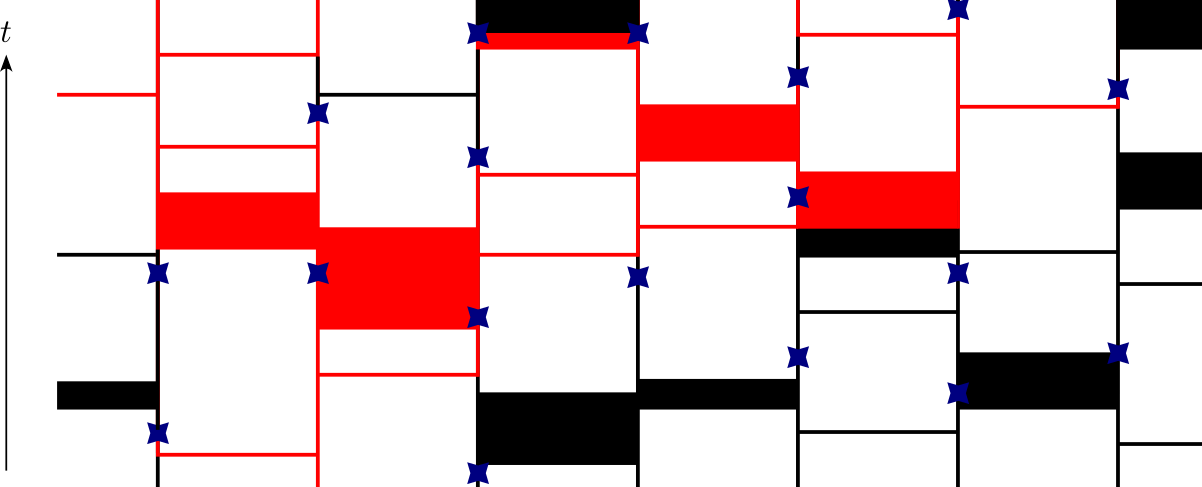}
    \caption{Graphical representation of a generalised contact process with undirected edges. Vertical lines represent vertices $v$ with crosses indicating recovery times $\Y_v$. The path of an infection is marked red and moves in the positive direction of time until it hits a recovery time, all while spreading to neighbouring vertices $v'$ along horizontal bars. The latter indicates infection times $\X_{(v,v')}$.}
    \label{fig:graphical-representation}
\end{figure}

\begin{definition}[Generalised contact process]\label{def:contact-process}
    Let $G=(V',E')$ be a graph with vertex set $V'$ and edge set $E'$. Let $\X,\Y\subset\R$ be random closed sets. We associate to vertices $v\in V'$ i.i.d.\ recovery times $\Y_v$ distributed like $\Y$ as well as to edges $e\in E'$ i.i.d.\ infection times $\X_e$ distributed like $\X$. Given an initial configuration of infected vertices $J_0\subset V'$, we define the set $J_t$ of infected vertices at time $t\geq0$ as follows:
    $v\in J_t$ if and only if there exists $v_0\in J_0$ and a finite path $\gamma=(v_0,v_1,\dots,v_k=v)$ together with contact times $0=t_{-1}\leq t_0\leq t_1,\dots,t_{k-1}\leq t_k=t$ such that,
    $$t_i\in \X_{(v_i,v_{i+1})}\quad \text{ for every } 0\leq i\leq k-1
    \qquad\text{and}\qquad
    [t_{i-1},t_i]\cap \Y_{v_i}=\emptyset \quad\text{ for every } 0\leq i \leq k.$$
    We say the contact process dies out if $J_t=\emptyset$ for some $t\geq0$.
\end{definition}
    Put differently, a vertex $v$ is infected at time $t$ if and only if there is a path of increasing infection times starting in some infected vertex $v_0$ and the corresponding vertices in the path do not recover in the time between being infected and passing on the infection. Depending on the setting, one may consider both directed and undirected graphs; all results in this section apply to both cases. An illustration is given in Figure~\ref{fig:graphical-representation}.
    
In the remainder, we explicitly assume $G=\Zz^d$ with the usual nearest-neighbour structure and the initial configuration $J_0=\{o\}\subset \Zz^d$ that only contains the origin.  We denote by $\PPP(\lambda)$ the homogeneous Poisson process of intensity $\lambda\ge 0$. Let us collect two general results before moving on to explicit examples. 

First, as presented in~\cite[Thm.~1]{fontes2020contactII}, let $\Y\subset\R$ be a renewal process with earliest recovery time after $t\geq0$ defined as
\begin{equation}\label{eq:earlist_recover}
 Z_t:=\min\{s\geq 0\colon s+t\in\Y\}.  
\end{equation}
If there exists $C>0$ such that $\E[Z_t]\leq C$ for every $t\geq 0$ and if $\X\sim \PPP(\lambda)$ with $\lambda < (2dC)^{-1}$, then the contact process on $\Z^d$ dies out almost surely.
The proof is based on a coupling to an i.i.d.~Bienaymé--Galton--Watson tree. In particular, the statement can be generalised to arbitrary random closed sets in which $Z_t$ is uniformly bounded.

Secondly, as seen in various settings~\cite{hilario2022contactprocess,fontes2023renewal,santos2023survival}, survival of the contact process in dimensions two or higher can be achieved in a simple way. The infection has enough space to sequentially discover new vertices on its path, thus eliminating dependencies. Our focus lies on the stationary case. The proof is given at the beginning of Section~\ref{sec:contact-processes}.

\begin{lemma}[Survival of infection in $d\geq 2$]\label{lem:survival-cp-2d}
    Let $p_{\mathsf{bond}},p_{\mathsf{site}}\in(0,1)$ be a pair of parameters for which the oriented Bernoulli bond-site percolation model on $\Ltwo$ contains an unbounded connected component. Let $\X$ and $\Y$ be stationary random closed sets with 
    \[
        \P(\X\neq\emptyset)>p_{\mathsf{bond}} \qquad\text{and}\qquad \P(0\notin\Y)>p_{\mathsf{site}}.
    \]
    Then, the corresponding generalised contact process on $\Zz^d,\,d\geq2$, with infection times $\X$ and recovery times $\mu^{-1}\Y$ survives with positive probability for all sufficiently small $\mu$. 
\end{lemma}

Now, let us move on to explicit examples. Consider a network of devices in which a virus spreads from device to device at random times. While the devices are patched on a regular schedule, thereby deleting the virus from the system, they still remain susceptible. This motivates considering a contact process with periodic recoveries. If $\Y=2\Z$ almost surely, then the infection will die out after two units of time have passed since all vertices recover simultaneously. Therefore, we study cases in which the vertices recover asynchronously.
\begin{thm}[Survival of contact processes with periodic recovery (CPPR)]\label{thm:CPPR-survival} 
Consider the generalised contact process on $\Zz$ with $\X\sim\PPP(\lambda)$ and either 
\begin{align}\label{eq:thm:CPPR-survival-continuous}\tag{Uni}
    \Y&=2 (\Z + U), \qquad \text{where $U$ is a uniform random variable on $[0,1)$, or}\\
    \Y&=2\Z + B, \qquad \text{where $B$ is a Bernoulli random variable with parameter $q\in(0,1)$.} \label{eq:thm:CPPR-survival-discrete}\tag{Ber}
\end{align}
Then, the contact process has a positive probability of survival for sufficiently large $\lambda$.
\end{thm}

We call this the {\em contact process with periodic recovery (CPPR)}. If $\lambda \leq (4d)^{-1}$, the infection dies out by~\cite[Thm.~1]{fontes2020contactII}, as described above. Establishing survival for large $\lambda$, which can be done by utilising Theorem~\ref{thm:OPRE-percolates}, then implies the existence of a non-trivial survival/extinction phase transition.

%Moving away from uniform shifts $U$ in~\eqref{eq:thm:CPPR-survival-continuous} to more general distributions on $[0,1)$ is feasible. As will become clear from the proof, what matters is the torus distance between neighbouring realisations.

\subsection{Contact processes in random spatial environment}
Finally, let us consider the setting of~\cite[Thm.~2]{bramson1991contact}, i.e.\ the contact process in random environment, which can be described as a generalised contact process, where $\X\sim\PPP(1)$ and $\Y$ is the following Cox point process: Let $p\in (0,1)$ and $\Delta,\delta>0$. Then, $\Y\sim\PPP(\delta)$ with probability \(p\) and $\Y\sim\PPP(\Delta)$ otherwise. Theorem~2 of~\cite{bramson1991contact} states that survival occurs for every $\delta<\delta_c(p,\Delta)$ for some critical $\delta_c(p,\Delta)>0$. We present an improved version of this statement in which the recovery rate is unbounded.

\begin{proposition}[Contact process in (spatial) random environment]\label{prop:CPRE-BDS}
     Let $\Delta>0$ be a random variable, with $\E[\Delta^{1+\eps}]<\infty$ for some $\eps>0$, and fix $\delta>0$ and $p\in(0,1)$. Let $(\Delta_x)_{x\in\Zz}$ be an i.i.d.\ family of random variables distributed according to $\mathcal{L}(\Delta)$, the law of $\Delta$. Consider the generalised contact process on $\Zz$ with $\X\sim\PPP(1)$ and recovery times  $(\Y_x)_{x\in\Zz}$ given by independent copies of the following Cox point processes:
    \begin{enumerate}[(i)]
        \item With probability $p$ (independent of everything else), let $\Y_x\sim\PPP(\delta)$.  
        \item Otherwise, let $\Y_x\sim\PPP(\Delta_x)$.
    \end{enumerate}
    Then, there exists some $\delta_c:=\delta_c(p,\mathcal L(\Delta))>0$ such that, for every $\delta<\delta_c$, the generalised contact process has a positive probability of survival. 
\end{proposition}
In essence, this is a generalised contact process with infection times $\X\sim\PPP(1)$ and recovery times $\Y$, where the Cox point process $\Y$ is given by a $\PPP$ of random intensity determined by $p$, $\delta$ and $\mathcal{L}(\Delta)$.
%In order to avoid confusion, the model is such that each vertex $x$ has its own independent $\Delta_x$ attached to its recovery time $\Y_x\sim\PPP(\Delta_x)$ with $\Delta_x$ following the law $\mathcal L(\Delta)$.

Let us also briefly sketch a second setting, which is a consequence of Proposition~\ref{prop:CPRE-BDS} after a time-rescaling.
\begin{cor}[Contact process in (spatial) random environment II]\label{cor:CPRE-2}
    Let $\Delta>0$ be a random variable, with $\E[\Delta^{1+\eps}]<\infty$ for some $\eps>0$. Let $(\Delta_x)_{x\in\Zz}$ be an i.i.d.\ family of random variables distributed according to $\mathcal{L}(\Delta)$. Consider the generalised contact process on $\Zz$ with $\X\sim\PPP(\lambda)$ and independent recovery times  $(\Y_x)_{x\in\Zz}$, where $\Y_x\sim\PPP(\Delta_x)$. Then, the generalised contact process has a positive probability of survival for every sufficiently large $\lambda$.
     %Consider the generalised contact process on $\N$ with $\X\sim\PPP(\lambda)$ and $\Y$ being a Cox point process defined as follows. Let $\Delta>0$ be a random variable with $\E[\Delta^{1+\eps}]<\infty$ for some $\eps>0$. Then, a realisation of $\Y$ is given by first taking a realisation of $\Delta$ and then $\Y\sim\PPP(\Delta)$. Then, the contact process is supercritical for $\lambda$ large enough.
\end{cor}
Again, this is a generalised contact process with infection times $\X\sim\PPP(\lambda)$ and recovery times $\Y$, where the Cox point process $\Y$ is given by a Poisson point process of random intensity $\Delta$. We omit the proof.
%%%%%%%%%%%%%%%%%%%%%%%%%%%%%%%%%%%%%%%%%%%%%%%%%%%%%%%%%%%%%%%%%%%%%%%%%%%%%%%%%%%%%%%%%%%%%%%
%%%%%         Proofs
%%%%%%%%%%%%%%%%%%%%%%%%%%%%%%%%%%%%%%%%%%%%%%%%%%%%%%%%%%%%%%%%%%%%%%%%%%%%%%%%%%%%%%%%%%%%%%%
%%%%%%%%%%%%%%%%%%%%%%%%%%%%%%%%%%%%%%%%%%%%%%%%%%%%%%%%%%%%%%%%%%%%%%%%%%%%%%%%%%%%%%%%%%%%%%%

\section{Proofs}
In this section we present all proofs. We start by showing the applications of Theorem~\ref{thm:OPRE-percolates} to the various contact processes that all rely on their own unique coupling. Afterwards, we present the proof of our main result on oriented percolation in Section~\ref{sec:OPRE}. 
\subsection{Survival of contact processes}\label{sec:contact-processes}
We start with the CPPR models, which feature some interesting behaviour due to their periodic recoveries, which we highlight next.
\begin{description}
    \item[Determinism:] Revealing a recovery time fixes it forever and, in particular, correlations do not decay over time.
    \item[Simultaneous healing:] Survival is difficult in blocks of vertices with similar recovery times as they all heal simultaneously. These blocks can be arbitrarily large.    
    \item[Time limitation:] The infection only has limited time to leave such a block before it simultaneously recovers. Hence, the probability of successfully crossing such a block of size $k$ scales roughly like $\lambda^k/k!$ in Model~\eqref{eq:thm:CPPR-survival-discrete}, in contrast to exponential decay in classical subcritical contact processes.
    %\item The CPPR is ergodic since it is mixing in diagonal directions -- just like many other percolation models with columnar disorder.
\end{description}

We start by proving Lemma~\ref{lem:survival-cp-2d}.

\begin{proof}[Proof of Lemma~\ref{lem:survival-cp-2d}] 
It suffices to show survival for the generalised contact process on the oriented north-east lattice on $\Zz^2$, which is equivalent to the $\Ltwo$ lattice. For this, we couple the contact process on $\Ltwo$ to oriented i.i.d.\ Bernoulli bond-site percolation also on $\Ltwo$. 
Since $\Y$ is almost-surely closed, we have
$$\bigcap_{k\in\N}\left\{ \Y\cap[0,1/k]\neq\emptyset\right\} = \left\{0\in\Y\right\}.$$
In particular, we have $\lim_{\eps\to0}\P(\Y\cap[0,\eps]=\emptyset)=\P(0\not\in\Y)>p_{\mathsf{site}}$. Now, fix $t>0$ large and \(\mu\) small such that
\[
	\P(\X\cap[0,t]\neq\emptyset)>p_{\mathsf{bond}}\qquad\text{and}\qquad \P(\Y\cap[0,2t\mu]=\emptyset)>p_{\mathsf{site}}.
\]
The coupling is constructed as follows. A vertex $(j,i)\in\Ltwo$ is called open if $\mu^{-1}\Y_{(j,i)}\cap[(j-1)t,(j+1)t]=\emptyset$. By the assumption above and stationarity, this happens with probability at least $p_{\mathsf{site}}$, independently of all other edges and vertices. Analogously, we call an edge $e$ linking $(j,i)$ to either $(j+1,i-1)$ or $(j+1,i+1)$ open if $\X_e\cap[jt,(j+1)t]\neq\emptyset$. By construction, any open path in the mixed bond-site percolation model corresponds to an infection path of the generalised contact process. The claim follows as, by the choice of $t$ and $\mu$, the associated bond-site percolation model contains an infinite connected component. 
\end{proof}

\subsubsection{CPPR -- continuous setting, Case~\eqref{eq:thm:CPPR-survival-continuous}}
\begin{figure}[ht]
    \centering
    \includegraphics[width=0.8\linewidth]{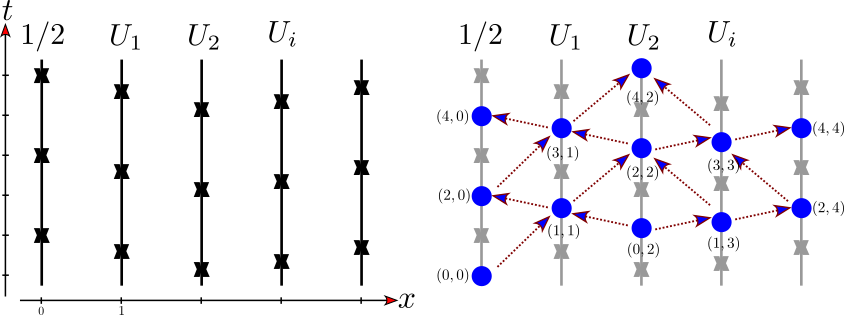}
    \caption{Left: Realisation of the random environment given by randomly shifted periodic recovery times $X_x=2(\Z+U_x)$. We assume without loss of generality that $U_0=1/2$. Right: Denoted in blue are the artificial lattice points required for the renormalisation scheme. Assuming $U_0=1/2$, the first lattice point is $(j,i)=(0,0)$ at location $(t,x)=(0,0)\in\R^2$. The other vertices are labelled with their $(j,i)$ coordinates.}
    \label{fig:CPPR-continuous-realisation}
\end{figure}

\begin{figure}[ht]
    \centering
    \includegraphics[width=0.8\linewidth]{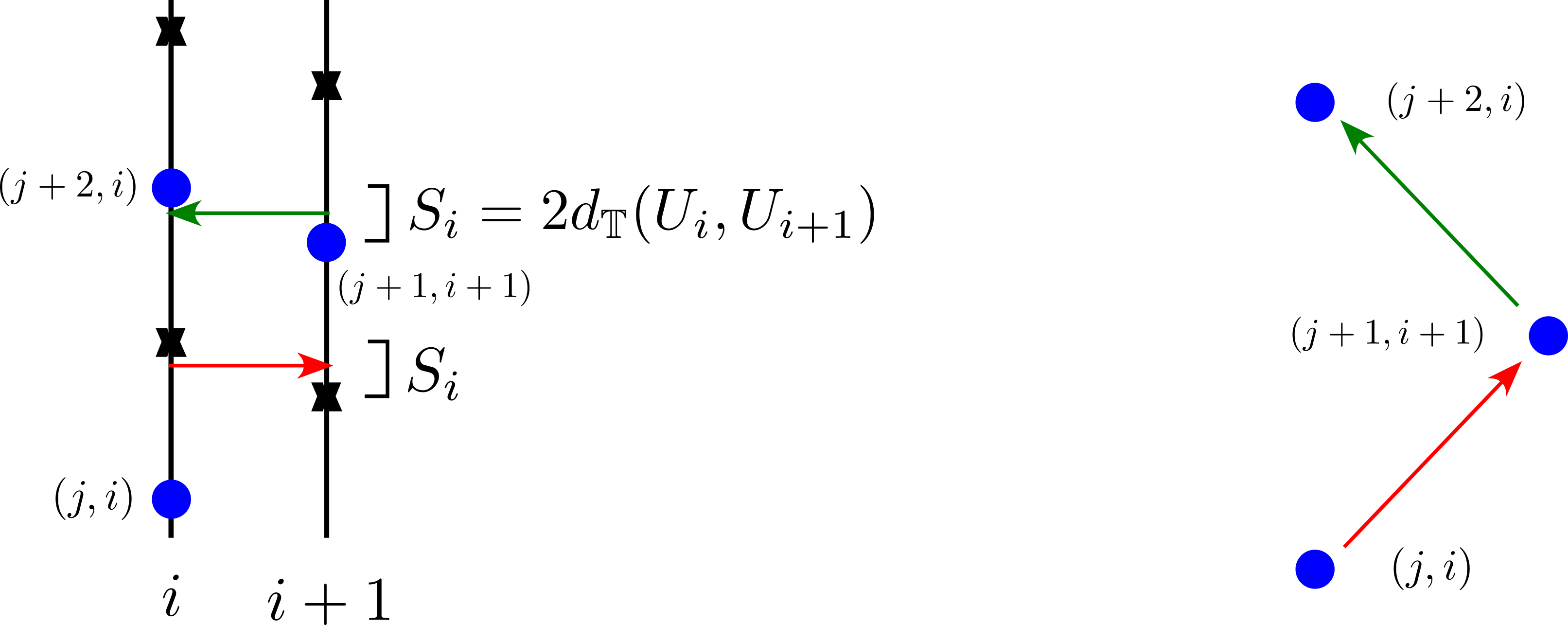}
    \caption{Left: Illustration of two infection events in the generalised contact process within their associated time intervals of length $S_i$. Right: Corresponding open edges in the coupled OPRE.}
    \label{fig:CPPR-continuous-edge-probability}
\end{figure}

\begin{proof}[Proof of Theorem~\ref{thm:CPPR-survival}, Case~\eqref{eq:thm:CPPR-survival-continuous}]
 Without loss of generality, we assume $U_0=1/2$, otherwise transform $U_i \mapsto U_i - U_0+1/2$ and start the process from a different time. Given a realisation of the $\Y_i$ (equivalently $U_i$), we assign vertices $(j,i)\in\Ltwo$ as in Figure~\ref{fig:CPPR-continuous-realisation}. Now we construct a coupling to the OPRE in which vertices are always considered open and edges are declared open as follows. Let $S_i:=2 d_\mathbb{T}(U_i,U_{i+1})$ be twice the torus distance between $U_i$ and $U_{i+1}$. Note that $S_i$ is still uniformly distributed on $[0,1)$. Vertex $(j,i)$ can infect vertex $(j+1,i+1)$ if an infection event occurs in the associated time interval of length $S_i$, see Figure~\ref{fig:CPPR-continuous-edge-probability}. If this is the case, we declare the corresponding edge $(j,i)\to(j+1,i+1)$ open. The same applies to the edge $(j+1,i+1)\to(j+2,i)$. Given $S_i$, these events both have probability $1-\exp(-\lambda S_i)$. Therefore, for a given environment $(S_i)_{i\in\Zz}$, the procedure yields independent bond percolation on $\Ltwo$.

For a uniform random variable $S$ on $[0,1)$, take spatial stretches 
\[
	\nu := - \log S.
\]
In particular, $\nu$ is an exponential random variable with parameter one and has exponential moments. 
Consider the connection function
\[
	\kappa_\lambda(s) := 1- \exp(-\lambda \e^{-s}),
\]
which is monotonically increasing in $\lambda$, and satisfies 
\[
	\kappa_\lambda(\nu) =  \kappa_\lambda(-\log S) = 1-\exp(-\lambda S),
\]
which is exactly the probability of an edge being open under the environment variable $S$. 
Hence, we must only verify Condition~\eqref{eq:OPRE-requirements-connection-kernel} of Theorem~\ref{thm:OPRE-percolates}. 
Clearly, $\kappa_\lambda$ is monotone and converges pointwise to $1$ as $\lambda\to\infty$. Next, we show the statement for $\lambda_0=2$ and sufficiently large $s$. This follows from the Taylor approximation 
\[
	\exp\left(-2 \e^{-s}\right) = 1 - 2 \e^{-s} \left(1 + h(s)\right),
\]
for some $h\colon\R\to\R$ with $h(s)\to0$ as $s\to\infty$. Take $s_0$ such that $h(s)\geq -1/2$ for every $s\geq s_0$. Then, for every $\lambda\geq2$, we have
\begin{equation}
    \begin{aligned}
        \kappa_\lambda(s) - \e^{-s} \geq \kappa_2(s) - \e^{-s} 
        = 1- \exp\left(-2 \e^{-s}\right) - \e^{-s} 
        =  \e^{-s} \left(2 - 1 + 2 h(s) \right) \geq 0.
    \end{aligned}
\end{equation}
Hence,~\eqref{eq:OPRE-requirements-connection-kernel} holds for all $s\geq s_0$. On the other hand, by choosing $\lambda$ sufficiently large, we also have $\kappa_\lambda(s_0)\geq \e^{-1}$ and in particular $\kappa_\lambda(s)\geq \e^{-s}$ for all $s\in [1,s_0]\cap\N$. 
Choosing $\nu_{i,i+1} =-\log(S_i)=-\log(2 d_\mathbb{T}(U_i,U_{i+1}))$ and $\xi_i\equiv 0$, provides a coupling of the CPPR model to an OPRE that satisfies the conditions of Theorem~\ref{thm:OPRE-percolates}. This finishes the proof.
\end{proof}

\subsubsection{CPPR -- discrete setting, Case~\eqref{eq:thm:CPPR-survival-discrete}}

\begin{figure}[ht]
    \centering
\includegraphics[width=0.8\linewidth]{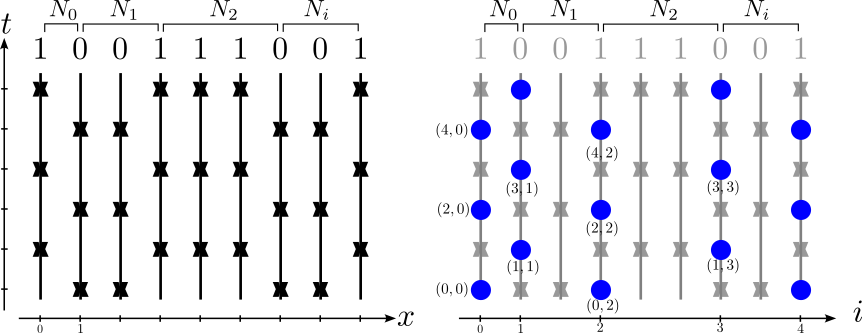}
    \caption{Left: Realisation of the random environment $(B_x)_{x\in\Zz}$, as illustrated as $1$s and $0$s, given by randomly shifted periodic recovery times. We consider intervals of consecutive $0$s ending in a $1$ and vice versa. Their lengths are denoted by $N_i$. These are independent geometric random variables with parameters alternating between $q$ and $1-q$. Right: We consider artificial lattice points (blue) in order to represent the renormalisation scheme. Those correspond exactly to vertices in the $\Ltwo$-lattice.}
    \label{fig:CPPR-discrete-realisation}
\end{figure}

\begin{figure}[ht]
    \centering
\includegraphics[width=0.65\linewidth]{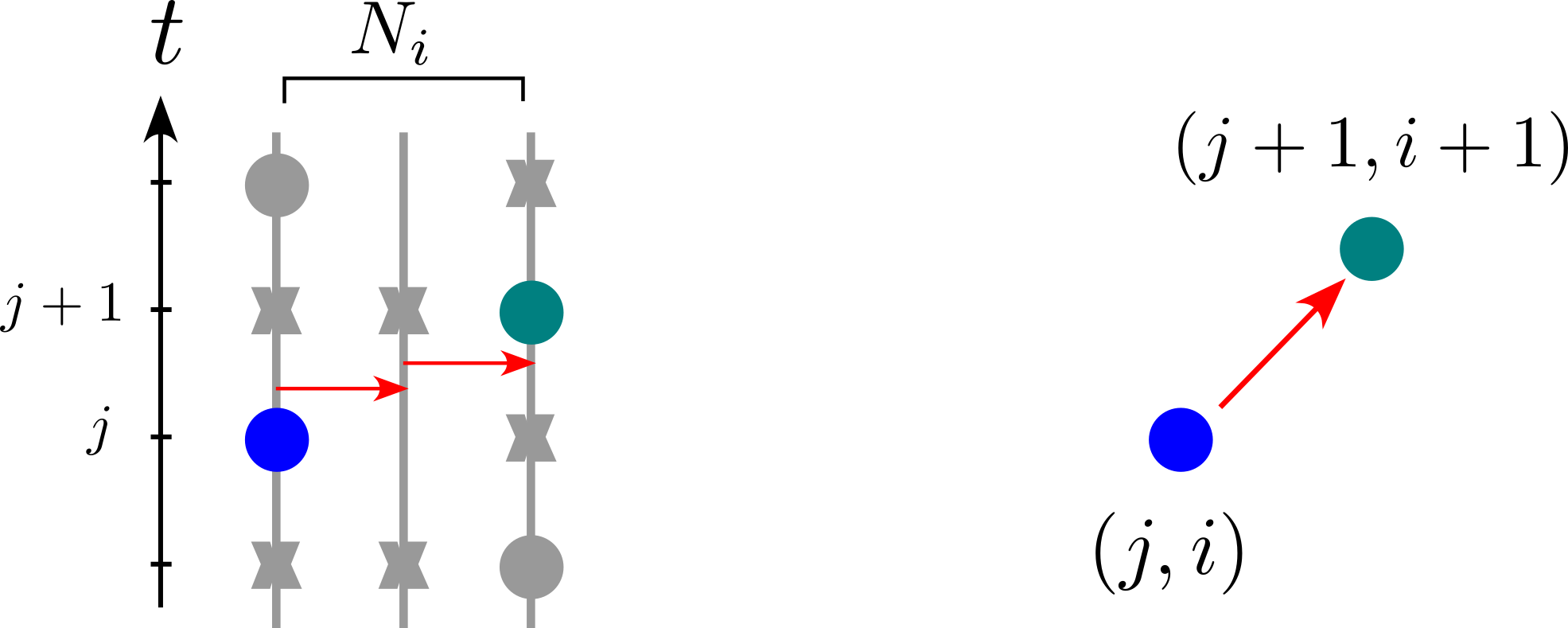}
    \caption{Left: Starting from the blue vertex $(j,i)$, an infection can reach the teal vertex $(j+1,i+1)$ if $N_i$ many infection arrows exist with the correct ordering. This is equivalent to the sum of $N_i$ i.i.d.\ exponential random variables of parameter $\lambda$ being less than $1$. This -- again -- is equivalent to having $N_i$ many Poisson points in an interval of length $1$. Right: An edge $(j,i)\to(j+1,i+1)$ is open in $\Ltwo$ if the prior event happens. This happens with probability $\P(A_{N_i})$.}
    \label{fig:CPPR-discrete-edge-probability}
\end{figure}

\begin{proof}[Proof of Theorem~\ref{thm:CPPR-survival}, Case~\eqref{eq:thm:CPPR-survival-discrete}]
First, we define the event $A_k$  of having $k$ exponential random variables with parameter $\lambda$ add up to at most $1$. This is equivalent to a Poisson point process of parameter $\lambda$ having at least $k$ points in an interval of length $1$. Hence,
\begin{equation}\label{eq:definition-A_k}
    \P(A_k)= 1-\e^{-\lambda}\sum_{i=0}^{k-1}\frac{\lambda^i}{i!} 
     =\e^{-\lambda}\sum_{i=k}^\infty \frac{\lambda^i}{i!}\geq \e^{-\lambda}\frac{\lambda^k}{k!}.
\end{equation}
Given $s\in\N$, we consider the connection function
\[
	\kappa_\lambda(s):= \P(A_{\lfloor\sqrt{s}\rfloor}) = 1-\e^{-\lambda}\sum_{i=0}^{\lfloor\sqrt{s}\rfloor-1}\frac{\lambda^i}{i!}.
\]
We assume $\lambda\geq100$, to facilitate calculation, and note that 
\[
	\kappa_\lambda(s)\geq 1/2 \geq \e^{-s} \quad \text{ for all } 1\leq s\leq \lambda^2,
\]
by the usual Poisson concentration inequalities. 
On the other hand, for $s\geq \lambda^2$ and writing $r=\lfloor\sqrt{s}\rfloor\geq \lambda-1$, Stirling's formula yields,
\begin{equation}
    \begin{aligned}
        \kappa_\lambda(s) 
        &
        	\geq \e^{-\lambda} \tfrac{\lambda^r}{r!}
        	\geq \tfrac{1}{2\sqrt{2\pi r}} \exp\left(-\lambda + r \log\lambda - r (\log r - 1) \right) %\quad\vert\,r\geq\lambda\\
        	\geq \tfrac{1}{10\sqrt{r}} \exp\left(- r (\log r - \log\lambda) \right) %\quad\vert\,r\geq\lambda\geq 100\\
       \\ &
       		\geq \exp\left(- r \log r \right) 
       		\geq \exp(-r^2) 
       		\geq \e^{-s},
    \end{aligned}
\end{equation}
where we used \(r+1\geq\lambda\geq 100\).

Now, let $K$ be a geometric random variable with $\P(K\geq \ell+1) = \max\{q,1-q\}^\ell$, and set $\nu:=K^2$, which has all polynomial moments. Then, $\P(A_K)=\kappa_\lambda(\nu)$ and we make use of Theorem~\ref{thm:OPRE-percolates} by coupling the CPPR in Case~\eqref{eq:thm:CPPR-survival-discrete} to the OPRE. First, we couple the CPPR to a bond percolation model on $\Ltwo$ as described in Figure~\ref{fig:CPPR-discrete-realisation} and~\ref{fig:CPPR-discrete-edge-probability}. 
By yet another coupling, we may assume that the $N_i$, as defined in Figure~\ref{fig:CPPR-discrete-realisation}, are independent geometric random variables with the same distribution as $K$ since that only increases the distances. However, this is exactly the OPRE with $\kappa_\lambda$ and $\nu$ as chosen above and vertices that are always open. This concludes the proof.
\end{proof}

\subsubsection{Contact processes in random environment}
\begin{proof}[Proof of Proposition~\ref{prop:CPRE-BDS}]
    The proof bears similarities to the proof of Theorem~\ref{thm:CPPR-survival}, Case~\eqref{eq:thm:CPPR-survival-discrete}. As before, we give a visual summary on the discretisation scheme in Figure~\ref{fig:discretisation-bds}. 
\begin{figure}[ht]
    \centering
\includegraphics[width=0.95\linewidth]{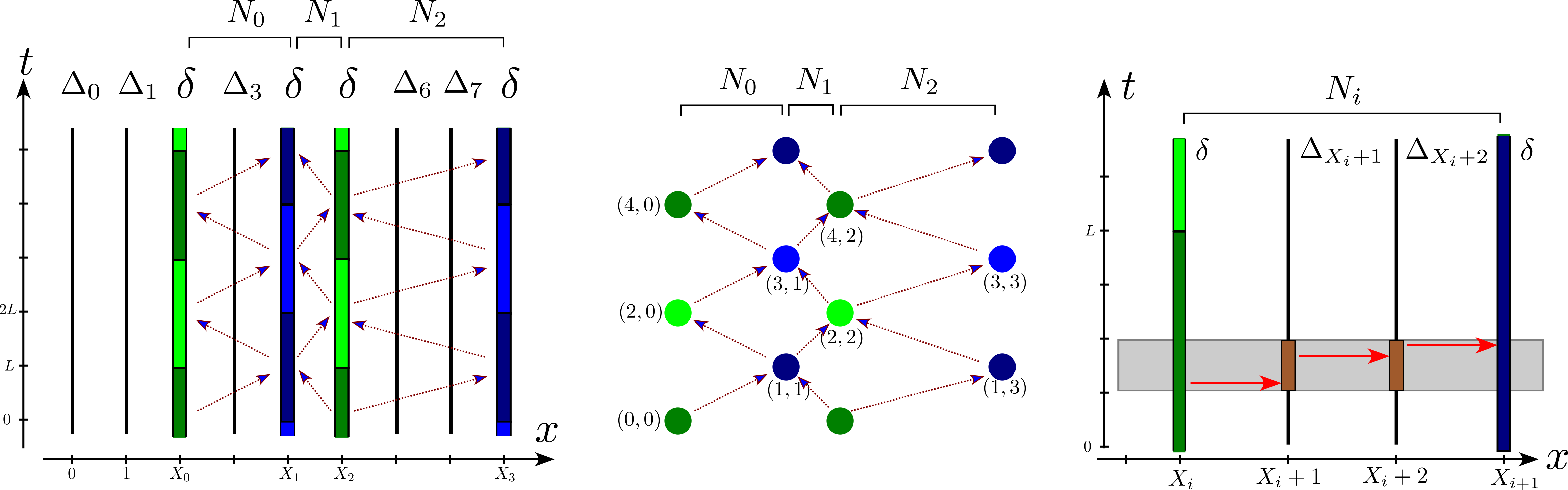}
    \caption{Left: We divide the process into temporal blocks of length $L$. Coloured boxes correspond to vertices in $\Ltwo$ with arrows inbetween indicating certain infection events. Middle: Coupling to a bond-site percolation model on $\Ltwo$. Vertices are open if the corresponding coloured box contains no recovery event. This happens with probability $\exp(-2L\delta)$. An edge $e$ is open if one of the $L$ many events $A_e(k)$ happens. Right:  Depicted is $A_e(1)$ for $e=(0,i)\to(1,i+1)$. The event happens if the red infection arrows exist in the proper order in the time interval $[1,2)$ and if no recovery events take place inside the brown boxes. 
    }
    \label{fig:discretisation-bds}
\end{figure}
    Let $(\Delta_x)_{x\in\Zz}$ be sequence of i.i.d.\ copies of $\Delta$ and let $\Y_x\sim\PPP(\Delta_x)$ with probability $1-p$, or $\Y_x\sim\PPP(\delta)$ otherwise. We declare a location \(x\) \emph{good} if its environment satisfies \(\Y_x\sim\PPP(\delta)\) and denote by \((X_i)_{i\in\Zz}\) the sequence of good locations in increasing order. 
   % We enumerate the locations of good environments in order as $(X_i)_{i\in\Zz}\subset\Zz$, that is, those with $\Y_{X_i}\sim\PPP(\delta)$. 
    We set $N_i=X_{i+1}-X_i$, \(i\in\Zz\) and note that these define a sequence of independent geometric random variables with parameter $p$.

Fix a length scale $L\in\N$, and discretise the model to an OPRE as follows. A vertex $(j,i)\in\Ltwo$ is called  open if $[(j-1)L,(j+1)L) \cap \Y_{X_i} = \emptyset$, which happens with probability
\[
	\P((j,i) \text{ is open})=\exp(-2\delta L).
\]
Note that these events are independent from each other for distinct vertices in $\Ltwo$. Given an edge $e=(j,b)\to(j+1,b+1)$ in $\Ltwo$ and $k\in\{0,\dots,L-1\}$, we define the event 
    \begin{equation}
    \begin{aligned}
        A_e(k):=\{ & \exists\, jL+k \leq t_1\leq \dots \leq t_{N_b} < jL+k+1 \colon \\
        &t_r \in \X_{X_b+r-1} \text{ and }\Y_{X_b+r}\cap [jL+k,jL+k+1)=\emptyset \, \text{ for all }  1\leq r\leq N_b-1 \}.
    \end{aligned}
    \end{equation}
    Loosely speaking, \(A_e(k)\) denotes the event that the vertex $(j,b)$ infects vertex $(j+1,b+1)$ with all infection events happen within the $k$-th time interval of length one. Given the random environment $(X_i)_{i\in\Zz},(\Delta_x)_{x\in\Zz}$, the existence of the necessary, ordered infections in the event $A_e(k)$ has probability $\e^{-1}\sum_{r\ge N_b}1/r!$, while the absence of recoveries has probability $\prod_{\ell=X_b+1}^{X_{b+1}-1}\e^{-\Delta_\ell}$. Thus, for all $k$, we have
    \[
        \P(A_e(k) \,\vert\, (X_i)_{i\in\Zz},(\Delta_x)_{x\in\Zz} )=\e^{-1}\sum_{r\ge N_b}\frac{1}{r!} \prod_{\ell=X_b+1}^{X_{b+1}-1}\e^{-\Delta_\ell}\,. 
    \]
    We declare the edge $e$ \emph{open} if $A_e(k)$ occurs for at least one $k\in\{0,\dots,L-1\}$. Conditioned on $(X_i)_{i\in\Zz}$ and $(\Delta_x)_{x\in\Zz}$, the events \(A_e(0),\dots,A_e(k)\) are i.i.d.\ for each edge \(e\), and the same holds true across different edges. Hence,
    \begin{equation}
        \begin{aligned}
            \P & \left(e \text{ is open} \,\vert\, (X_i)_{i\in\Zz},(\Delta_x)_{x\in\Zz}\right) 
            = 
            	1-\left[1-\P(A_e(0) \,\vert\, (X_i)_{i\in\Zz},(\Delta_x)_{x\in\Zz})\right]^L \\
     			&= 1 - \Big(1-   \e^{-1}\sum_{r\ge N_b}\frac{1}{r!} \cdot \prod_{\ell=X_b+1}^{X_{b+1}-1}\e^{-\Delta_\ell}  \Big)^L
     			\geq 1 - \Big(1-   \exp\Big(- N_b\log N_b - \sum_{\ell=X_b+1}^{X_{b+1}-1}\Delta_\ell\Big)  \Big)^L.
        \end{aligned}
    \end{equation}
    %Here, the part in the probability of $A_e(0)$ that deals with the process $\X$ is represented by $\e^{-1}\sum_{r\ge N_b} 1/r!$ and the part that deals with $\Y$ is represented by $\prod_{\ell=X_b+1,\dots,X_{b+1}-1}\e^{-\Delta_\ell}$.
    Choosing $\delta=L^{-2}$, this provides a coupling of the generalised contact process to an OPRE, where the edges use the connection function $\kappa_L(s)=1-(1-\exp(-s))^L$ with stretches $\nu=N_0\log N_0 + \sum_{\ell=0,\dots, N_0}\Delta_\ell$, and the vertices use the constant connection function $\kappa_L=\exp(-2L^{-1})$.
    Since
    \[
    	\E[\nu^{1+\eps}] 
			\leq 2^{1+\eps}\Big(\E[(N_0\log N_0)^{1+\eps}]+\E[N_0^{1+\eps}\sum^{N_0}_{r=0}\Delta_r^{1+\eps}] \Big)
            \leq 2^{1+\eps}\left(\E[N_0^{2+2\eps}]+\E[N_0^{2+\eps}]\E[\Delta_0^{1+\eps}] \right)
    		%\leq \E[N_0^{2+\eps}]\E[\Delta_0^{1+\eps}] 
    		< \infty,
   	\]
    Theorem~\ref{thm:OPRE-percolates} applies. This shows survival for sufficiently large $L$ (resp.\ small $\delta$) as each infinite oriented path starting from the origin in the respective OPRE corresponds to an infection path in the contact process.
   \end{proof}

\subsection{Oriented percolation in random environment (OPRE)}\label{sec:OPRE}
This section is dedicated to the proof of Theorem~\ref{thm:OPRE-percolates}. It is based on the static multiscale renormalisation scheme presented in~\cite{hilario2023phasetransition}, modified to account for the oriented bond-site percolation case. As such, Section~\ref{subsec:OPRE-renewal} introduces three main components: a stationary spatial embedding of the OPRE into $\Zz^2$, the notion of good (bad) blocks at different scales, as well as the estimate on the probability of a block being good. 
The alternative description of the OPRE is needed to enable a static renormalisation scheme, i.e., the blocks are \emph{not} defined via the random environment, only their state of being good/bad.
Section~\ref{subsec:OPRE-crossings} deals with estimates on long horizontal and vertical rectangle crossings (Lemmas~\ref{lem:OPRE-induction-horizontal-crossing} and~\ref{lem:OPRE-induction-vertical-crossing}).  Due to planarity, we can then construct an infinite, oriented open path. This proves the main Theorem~\ref{thm:OPRE-percolates}.

For those familiar with~\cite{hilario2023phasetransition}, we want to briefly highlight the main differences. Those are best described by comparing the figures in the current paper with~\cite{hilario2023phasetransition}. The first adjustment lies in the alternative stationary embedding of the graph (Figure~\ref{fig:OPRE-alternative}) due to the mixed bond-site percolation setting. More precisely, the stretches in the vertices also have an impact on the embedding and we attach those stretches to the right of the vertices. As a side effect of this change, all crossing events will disregard the right-most column of available vertices (Figure~\ref{fig:rectangle-crossings}). The oriented percolation setting impacts some arguments for vertical crossings (Figure~\ref{fig:vertical_crossing_renorm} vs.\ \cite[Fig.~7]{hilario2023phasetransition}) as well as both the construction of horizontal crossings (Figure~\ref{fig:LRC} vs.\ \cite[Fig.~7]{hilario2023phasetransition}) and the construction of the infinite open path (Figure~\ref{fig:construction_infinite_cluster} vs.\ \cite[Fig.~8]{hilario2023phasetransition}). The calculations remain mostly intact. Notably, the strengthened statement in terms of connection functions in Theorem~\ref{thm:OPRE-percolates} readily follows from the unaltered proof.

\subsubsection{Renewal process, alternative embedding, vertical scales}\label{subsec:OPRE-renewal}
\begin{figure}[ht]
    \centering
    \includegraphics[width=0.85\linewidth]{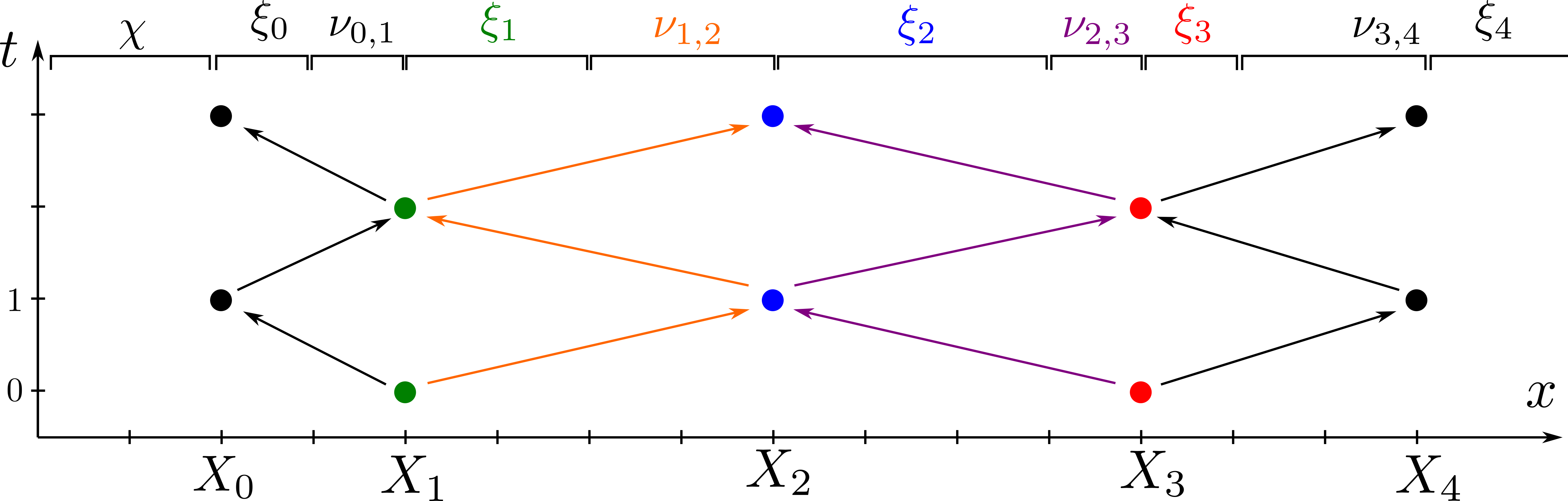}
    \caption{Stationary embedding of the graph in $\Zz^2$. Depicted is a realisation with $B=1$.}
    \label{fig:OPRE-alternative}
\end{figure}

First of all, we observe that we may set $\sigma=1$ in Condition~\eqref{eq:OPRE-requirements-connection-kernel} simply by replacing $\xi$ with $\sigma\xi$ (and $\nu$ with $\sigma\nu$), which we will assume from now on in the proofs. Without loss of generality, we assume \(\xi_i,\nu_{i,i+1}\in\N\) as we may replace \(\xi_i\) and \(\nu_{i,i+1}\) by $\lceil \xi_i \rceil$ and \(\lceil \nu_{i,i+1}\rceil\), respectively, without violating any assumption. Consider the stationarised renewal process with interarrival distribution $\xi_0+\nu_{0,1}$, i.e.
\begin{equation}\label{eq:OPRE-stationarised-renewal-X}
    X=\{X_i\}_{i\in\N} = \Big\{\chi + \sum_{k=0}^{i-1} \left(\xi_k + \nu_{k,k+1} \right)\Big\}_{i\in\N},
\end{equation}
where \(\chi\) is the difference of the forward-waiting time and \(\xi_0+\nu_{0,1}\).   
Let us slightly adapt the lattice in Definition~\ref{def:OPRE-model} in order to employ the static multiscale-renormalisation scheme of~\cite{hilario2023phasetransition}. More precisely, we introduce an additional coin flip that randomises the time of appearance of the first, left-most point of the lattice to be either \(t=0\) or \(t=1\). This contrasts the original definition, where said vertex always appears at time \(t=0\), and stationarises the construction, therefore leading to the following definition, see Figure~\ref{fig:OPRE-alternative}.

\begin{definition}[Alternative stationary spatial embedding]\label{def:OPRE-alternative-embedding} 
Let $B\in\{0,1\}$ be a Bernoulli random variable with parameter $1/2$ independent of everything else. Given (a realisation of) $B$ and $X$, consider the graph $G=(V,E)$ with vertex set
    \[
    	V=\big\{(t,X_i)\in\Zz^2 \colon t\in\Zz, t+i+B \text{ even}\big\}
    \]
and edge set
    \[
    	E=\big\{(t,X_i)\to(t+1,X_{i\pm 1}) \colon  (t,X_i)\in V\big\}.
    \]
Consider $\kappa_\lambda\colon[0,\infty)\to[0,1]$ as in Definition~\ref{def:OPRE-model}.
We consider the independent bond-site percolation model on $G$, where for either $(t,X_{i})\in V$ or $(t,X_{i+1})\in V$, the edges \((t,X_{i})\to(t+1,X_{i+1})\) and \((t,X_{i+1})\to(t+1,X_{i})\), respectively, are open with probability \(\kappa_\lambda(\nu_{i,i+1})\).
    %\[
    %	\begin{split}
    %		&
    %			\P_\lambda((t,X_{i})\to(t+1,X_{i+1}) \text{ is open}) := \kappa_\lambda(\nu_{i,i+1}) \qquad\text{and }
    %		\\ &
    %			\P_\lambda((t,X_{i+1})\to(t+1,X_{i}) \text{ is open}) := \kappa_\lambda(\nu_{i,i+1}),\qquad\text{  respectively}
    %	\end{split}
    %\]
    Similarly, vertices $(t,X_i)\in V$ are open independently with probability \(\kappa_\lambda(\xi_i)\).
    %\[
    %	\P_\lambda((t,X_i) \text{ is open}) := \kappa_\lambda(\xi_i).
    %\]
    We denote by $\P_\lambda^X(A):=\P_\lambda(A|X)$ the probability measure, conditionally given $X$ (as defined in~\eqref{eq:OPRE-stationarised-renewal-X}). 
\end{definition}

We see via coupling that almost-sure percolation of the OPRE in Definition~\ref{def:OPRE-model} implies almost-sure percolation of the model in Definition~\ref{def:OPRE-alternative-embedding} and vice versa. Hence, it suffices to show percolation in the latter model for almost-every realisation of $X$.

A key step is to identify good and bad spatial columns of the graph at different scales, for which we require the following quantities.

%%%%%%%%%%%%%%%%%%%%%%%%%%%%%%%%%%%%%%%%%%%%%%%%%%%%%%%%%%%%%%%%%%%%%%%%%%%%%%
\begin{parameters}\label{par:one} ~\
\begin{enumerate}[(i)]
    	\item 
    		$\eps>0$ is chosen sufficiently small such that $\E[\xi_0^{1+\eps}],\E[\nu_{0,1}^{1+\eps}]<\infty$, i.e., \eqref{eq:OPRE-stretches-moment-condition} is fulfilled.
    	\item 
    		Take any $\alpha\in(0,\eps/2]$. This parameter governs the probability of bad blocks.
    	\item 
    		Take any $\gamma\in(1,1+\alpha/(\alpha+2)]$. This parameter governs the rate of growth of scales.
    	\item 
    		Take $L_0\in\N$ such that $L_0^{\gamma-1}\geq5$ and Lemma~\ref{lem:OPRE-prob-bad-interval} below is satisfied. $L_0$ serves as the large initial block size.
    	\item 
    		We define $L_k:=L_{k-1}\lfloor L_{k-1}^{\gamma-1}\rfloor$. Notably, $L_k$ grows superexponentially and \begin{equation}\label{eq:estimate-L_k}
        		(1/2)^k L_0^{(\gamma^k)}\leq (1/2) L^\gamma_{k-1} \leq L_k \leq L_{k-1}^\gamma \leq L_0^{(\gamma^k)}.
    		\end{equation}
\end{enumerate}
\end{parameters}
%%%%%%%%%%%%%%%%%%%%%%%%%%%%%%%%%%%%%%%%%%%%%%%%%%%%%%%%%%%%%%%%%%%%%%%%%%%%%%

We start by defining \emph{good} blocks. As usual, each block that is not good is called \emph{bad}.

\begin{definition}[Good blocks] At fixed scale $k\in\N$, we partition $\Zz$ into blocks of length $L_k$, that is,
    \[
    	I_{k,i} := [iL_k, (i+1)L_k), \quad i\in\Zz,
    \]
    so that each scale-\((k+1)\) block $I_{k+1,i}$ consist of $\lfloor L_{k}^{\gamma-1}\rfloor$ many scale-$k$ blocks, i.e., $\#\{m\in\Z_{\geq0} \colon  I_{k,m}\cap I_{k+1,i}\neq\emptyset\} = \lfloor L_{k}^{\gamma-1}\rfloor$.
    On the initial scale, we call a block $I_{0,i}$ {\em good} if $I_{0,i}\cap X \neq \emptyset$. For general $k\in\N$, we call the scale-\(k\) block $I_{k,i}$ {\em good} if it either
    \begin{enumerate}[(i)]
        \item contains no more than one bad scale-\((k-1)\) sub-block or      
        \item if it contains exactly two consecutive bad sub-blocks, i.e., there exists \(m\) such that \(I_{k-1,m}, I_{k-1,m+1}\subset I_{k,i}\) are bad but all other sub-blocks \(I_{k-1,n}\subset I_{k,i}\) are good themselves.   
    \end{enumerate}
\end{definition}

The natural first step is to determine the probability of a given block to be bad, which can be found in~\cite[Lem.~3.1]{hilario2023phasetransition}. The result relies on a decoupling inequality, which is a consequence of the moment condition~\eqref{eq:OPRE-stretches-moment-condition} and is also the only situation where~\eqref{eq:OPRE-stretches-moment-condition} is used. We omit the proof as it is identical to the one in \cite{hilario2023phasetransition}.

\begin{lemma}[Probability of bad blocks at scale $k$, {\cite[Lem.~3.1]{hilario2023phasetransition}}]\label{lem:OPRE-prob-bad-interval}
    Let $\eps,\alpha,\gamma$ be as in Parameters~\ref{par:one}. Then, one can choose $L_0=L_0(\eps,\alpha,\gamma)$ sufficiently large such that the following holds. For every $k\in\N$ and $i\in\Zz$, we have that
    $$\P(I_{k,i}\text{ is bad}) \leq L_k^{-\alpha}.$$
\end{lemma}

% \begin{remark}[Regarding site percolation]\label{rem:OPRE-site-percolation}
%     \begin{figure}[ht]
%         \centering
%         \includegraphics[width=1.0\linewidth]{img/OPRE-site-percolation-embedding.png}
%         \caption{In the case of site percolation, an different embedding is needed as shown in the right picture. $\chi$ and $\xi_1$ also have to be adjusted accordingly so that the resulting picture is shift-invariant in space.}
%         \label{fig:OPRE-site-percolation}
%     \end{figure}
%     As previously mentioned, the same proof framework also allows to show survival of the infection when the stretches are placed on the vertices rather than edges of the graph. A slightly different embedding is needed as seen depicted in Figure~\ref{fig:OPRE-site-percolation}. Here, the distance from one site $i$ to $i+1$ is equal to the sum of their stretches, i.e.~$\xi_i+\xi_{i+1}$. To obtain a stationary picture, $\chi$ and $\xi_1$ also have to be adjusted accordingly. In turn, the proof of Lemma~\ref{lem:OPRE-prob-bad-interval} also needs some slight modifications in this setting.
% \end{remark}

%<<<<<<< Updated upstream
%\subsubsection{Horizontal scale, crossing events}\label{subsec:OPRE-crossings}
%Next, we study the probability of rectangle crossings in the oriented case. Given a rectangle $R:=%[t_0,t_1]\times[a,b]$, with $t_0<t_1$ and $a<b$, we will consider crossings in the subgraph of $G$ spanned by the vertices $V\cap R$. This graph does not necessarily have vertices in $\Zz\times \{a\}$ nor $\Zz\times \{b\}$, see Figure~\ref{fig:rectangle-crossings}.
%=======
\subsubsection{Horizontal scale, crossing events}\label{subsec:OPRE-crossings}
\begin{figure}
    \centering
    \includegraphics[width=1.0\linewidth]{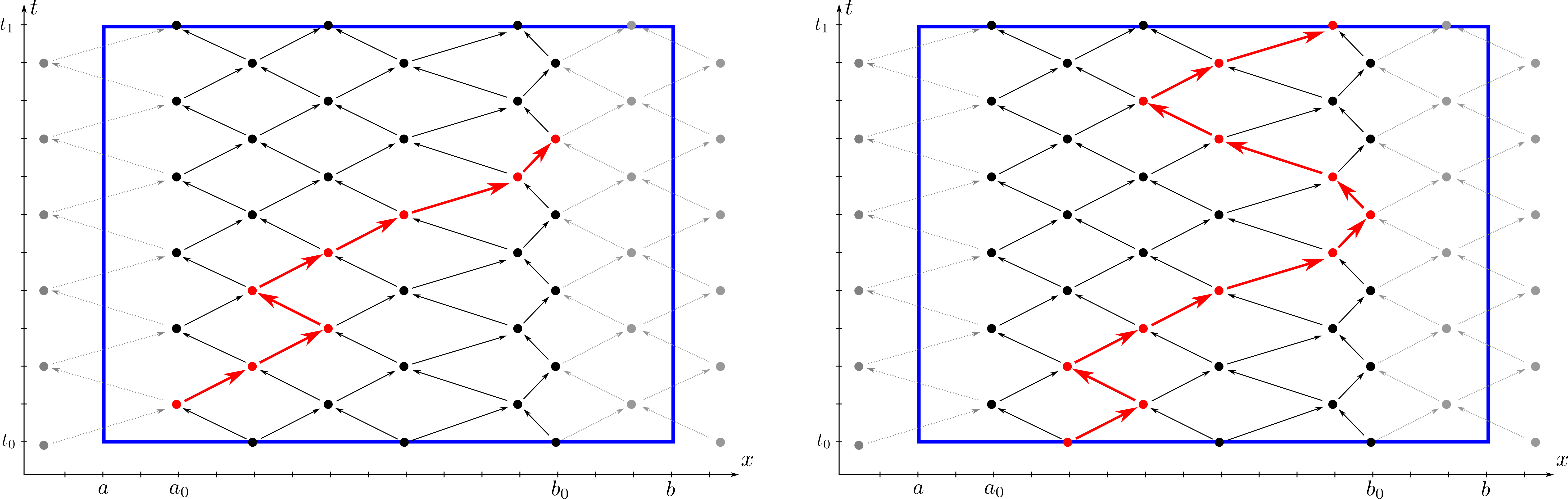}
    \caption{Depictions of the rectangle $R=[t_0,t_1]\times[a,b]$ (blue) together with the events $\LRC(R)$ and $\DTC(R)$. The corresponding open paths are marked in red.}
    \label{fig:rectangle-crossings}
\end{figure}

In this section, we study the probability of rectangle crossings in the oriented case. Given a rectangle $R:=[t_0,t_1]\times[a,b]$, %with $t_0<t_1$ and $a<b$, 
we consider crossings in the induced subgraph of $G$. %spanned by the vertices $V\cap R$. 
Note that this graph does not necessarily have vertices in $\Zz\times \{a\}$ nor $\Zz\times \{b\}$, see Figure~\ref{fig:rectangle-crossings}.
%>>>>>>> Stashed changes

Given such a rectangle $R$, we denote the left and right-most valid coordinates by $a_0$, respectively $b_0$. That is,
\begin{equation}\label{eq:def:reduced-rectangle}
	a_0 := \min\{X_i \colon X_i \geq a\},
	\qquad
	b_0 := \max\{X_{i-1} \colon X_i \leq b\},
	\qquad\text{and}\qquad
	R_0 := [t_0,t_1]\times [a_0,b_0].    
\end{equation}
Note that we do \emph{not} use the right-most column of vertices. Also, we will restrict to \emph{good} rectangles at scales $k\geq1$, guaranteeing the existence of at least three columns of vertices and, in particular, $a\leq a_0<b_0<b$. Hence, there is always a positive probability of traversing vertically. Furthermore, we choose a rectangle's height to be larger than its width, so that horizontal crossings have positive probability as well. More formally, we define
\begin{align*}
	\LRC(R)
	&
		:= \{\exists \text{ open path }[t_0,t_1]\times \{a_0\}\leadsto [t_0,t_1]\times \{b_0\} \text{ inside } R_0\},
	\\
	\RLC(R)
	&
		:= \{\exists \text{ open path }[t_0,t_1]\times \{b_0\}\leadsto [t_0,t_1]\times \{a_0\}  \text{ inside } R_0 \},\text{ and}
	\\
	\DTC(R)
	&
		:= \{\exists \text{ open path }\{t_0\}\times[a_0,b_0]\leadsto \{t_1\}\times[a_0,b_0]  \text{ inside } R_0 \}.
\end{align*}

Let us introduce next a scale-\(k\) rectangle and explain how it relates to the previously defined blocks. As supercritical contact processes typically survive exponentially long in the underlying graph size~\cite{schapira_extinction_2017,jahnel2025phasetransitionscontactprocesses}, we choose the height (time coordinate) of a rectangle exponentially in its width. To this end, we introduce another set of parameters, extending those of Parameters~\ref{par:one}.

\begin{parameters}\label{par:two} ~\
\begin{enumerate}[(i)]
    \item 
    	Take any $\mu\in \left(\gamma^{-1},1\right)$ to govern the height of scale-\(k\) rectangles.    
    \item 
    	Take any $\beta\in(\gamma\mu-\gamma+1,1)$, which governs the probability of crossings in good rectangles. Note that    \begin{equation} \label{eq:OPRE-beta-property}
        		\beta+\gamma-1>\max\{\gamma\beta,\gamma\mu\}.
    		\end{equation}
    \item Take any $H_0 = L_0$ as an initial height.
    \item Set $H_k:=2 \lfloor L_{k-1}^{\gamma-1} \rfloor \lceil \exp(L_k^\mu) \rceil H_{k-1}$, to determine the height of a scale-\(k\) rectangle.
\end{enumerate}
\end{parameters}

For horizontal crossings at scale $k\in\N$, we consider rectangles of the form
\begin{equation}\label{eq:OPRE-def-hor-rectangle}
    R^{\hor}_k(j,i) := [jH_k,(j+1)H_k]\times [iL_k,(i+2)L_k]=[jH_k,(j+1)H_k] \times (I_{k,i}\cup I_{k,i+1}),
\end{equation}
while, for vertical crossings, we use
\begin{equation}\label{eq:def-ver-rectangle}
    R^{\ver}_k(j,i) := [jH_k,(j+2)H_k]\times [iL_k,(i+1)L_k] = [jH_k,(j+2)H_k] \times I_{k,i}.
\end{equation}

Let $X$ be some realisation of the environment. We aim to appropriately bound the probability of these crossings happening in good intervals, given the environment. To this end, recall that we write \(\P^X_\lambda\) for the conditional probability measure, given the environment \(X\), and define
\[
	q_k(\lambda,j,i):= \max\Big\{   \max_{\substack{X\colon I_{k,i}, I_{k,i+1}\text{ good}}} \P^X_\lambda\big(\neg LRC(R^{\hor}_k(j,i))\big), \, \max_{\substack{X\colon I_{k,i}\text{ good}}} \P^X_\lambda\big(\neg DTC(R^{\ver}_k(j,i))\big) \Big\}.
\]
Note that the same bound holds for right-left-crossings. As $H_k$ and $L_k$ are always even, stationarity yields
$q_k(\lambda):=q_k(\lambda,0,0)=q_k(\lambda,j,i)$ for all $j,i\in\Zz$. In the following lemma, which is a modification of~\cite[Lem.~3.2]{hilario2023phasetransition}, we derive the crucial bounds on \(q_{k}\).

\begin{lemma} \label{lem:OPRE-main-induction}
    For every sufficiently large $\lambda$ depending on $\gamma,\, L_0,\, \mu,\, \beta$, we have
\begin{equation}\label{eq:OPRE-multiscale-prob-estimate}
          q_k(\lambda)\leq\exp(-L_k^\beta) \quad \forall k\geq1.
     \end{equation}
    % There exists $c_3=c_3(\gamma,L_0,\mu,\beta)\in\N$ such that, for all $k\geq c_3$, the following holds. If 
    % \begin{equation}\label{eq:OPRE-multiscale-prob-estimate}
    %     q_k(\lambda)\leq\exp(-L_k^\beta),
    % \end{equation}
    % then also
    % \[
    % 	q_{k+1}(\lambda)\leq\exp(-L_{k+1}^\beta).
    % \]
\end{lemma}

The proof is a direct consequence of the following two lemmas, which are analogous to~\cite[Lem.~3.3, 3.4]{hilario2023phasetransition}. %The proof proceeds inductively by considering the relevant horizontal and vertical crossings. As such, it follows immediately from the following two lemmas by taking $c_3:=\max\{c_4,c_5\}$ from their corresponding statements. Due to~\eqref{eq:OPRE-requirements-connection-kernel}, we can guarantee~\eqref{eq:OPRE-multiscale-prob-estimate} for $1\leq k\leq c_3$ by choosing $\lambda$ sufficiently large.

\begin{lemma}[Probability of horizontal crossings]\label{lem:OPRE-induction-horizontal-crossing}
    There exists $c_4=c_4(\gamma,L_0,\mu,\beta)\in\N$ such that, for all $k\geq c_4$, we have that
    \[
    	q_k(\lambda)\leq\exp\big(-L_k^\beta\big)\quad \Longrightarrow \quad
    \max_{\substack{X\colon I_{k,i},I_{k,i+1}\text{ good}}} \P^X_\lambda\big(\neg \LRC(R^{\hor}_{k+1}(0,0))\big)  \leq \exp\big(-L_{k+1}^\beta\big).
    \]
    The same holds true for right-left crossings.
\end{lemma}

\begin{lemma}[Probability of vertical crossings]\label{lem:OPRE-induction-vertical-crossing}
    There exists $c_5=c_5(\gamma,L_0,\mu,\beta)\in\N$ such that, for all $k\geq c_5$, we have that
    \[
    	q_k(\lambda)\leq\exp\big(-L_k^\beta\big)\quad \Longrightarrow \quad
    \max_{\substack{X\colon I_{k,i}\text{ good}}} \P^X_\lambda\big(\neg \DTC(R^{\ver}_{k+1}(0,0))\big)  \leq \exp\big(-L_{k+1}^\beta\big).
   	\]
\end{lemma}

The proofs are given in the next section. Before doing so, we explain how these lemmas imply Lemma~\ref{lem:OPRE-main-induction}.
\begin{proof}[Proof of Lemma~\ref{lem:OPRE-main-induction}]
First note that there are only finitely many configurations of $X$ inside a fixed rectangle at scale $k$. Using~\eqref{eq:OPRE-requirements-connection-kernel}, we further have $q_k(\lambda)\to0$ as $\lambda\to\infty$ for any fixed $k\geq1$. The proof now finishes by induction. First, set \(c_3=\max\{c_4,\, c_5\}\) and choose \(\lambda\) sufficiently large to guarantee \(q_k(\lambda)<\exp(-L_k^{\beta})\) for all \(k\leq c_3\). For $k\geq c_3$, the claim follows from inductively using Lemma \ref{lem:OPRE-induction-horizontal-crossing} and \ref{lem:OPRE-induction-vertical-crossing}. 
\end{proof}

\subsubsection{Probability of crossing events}
%Let us prove the Lemmas~\ref{lem:OPRE-induction-horizontal-crossing} and~\ref{lem:OPRE-induction-vertical-crossing} now. 
\begin{proof}[Proof of Lemma~\ref{lem:OPRE-induction-horizontal-crossing}] Fix an environment $X$ where the blocks $I_{k+1,0},I_{k+1,1}$ are good. Divide the rectangle $R^\hor_{k+1}(0,0)$ into disjoint strips of height $2\lfloor L_k^{\gamma-1} \rfloor H_{k}$ and try to traverse each of them from left to right. By definition of $H_{k+1}$, there are $\lceil \exp(L_{k+1}^\mu)\rceil$ many of those strips and hence that many independent trials to obtain a crossing. More precisely, consider the event.
\[
	G_{k+1}:=\LRC([0,2\lfloor L_k^{\gamma-1}\rfloor H_{k}] \times [0,2L_{k+1}]). 
\]
In good blocks, we use bottom-top crossings and planarity of the underlying graph to connect left-right crossings. Bad areas are attempted to be traversed as fast as possible. Recall that, at stage \(k\), there are no more than two bad areas consisting of at most two consecutive bad blocks while there are $2\lfloor L_k^{\gamma-1} \rfloor$-many scale-\(k\) blocks in total. The scheme is depicted in Figure~\ref{fig:LRC}.

\begin{figure}[ht]
    \centering
    \includegraphics[width=0.9\linewidth]{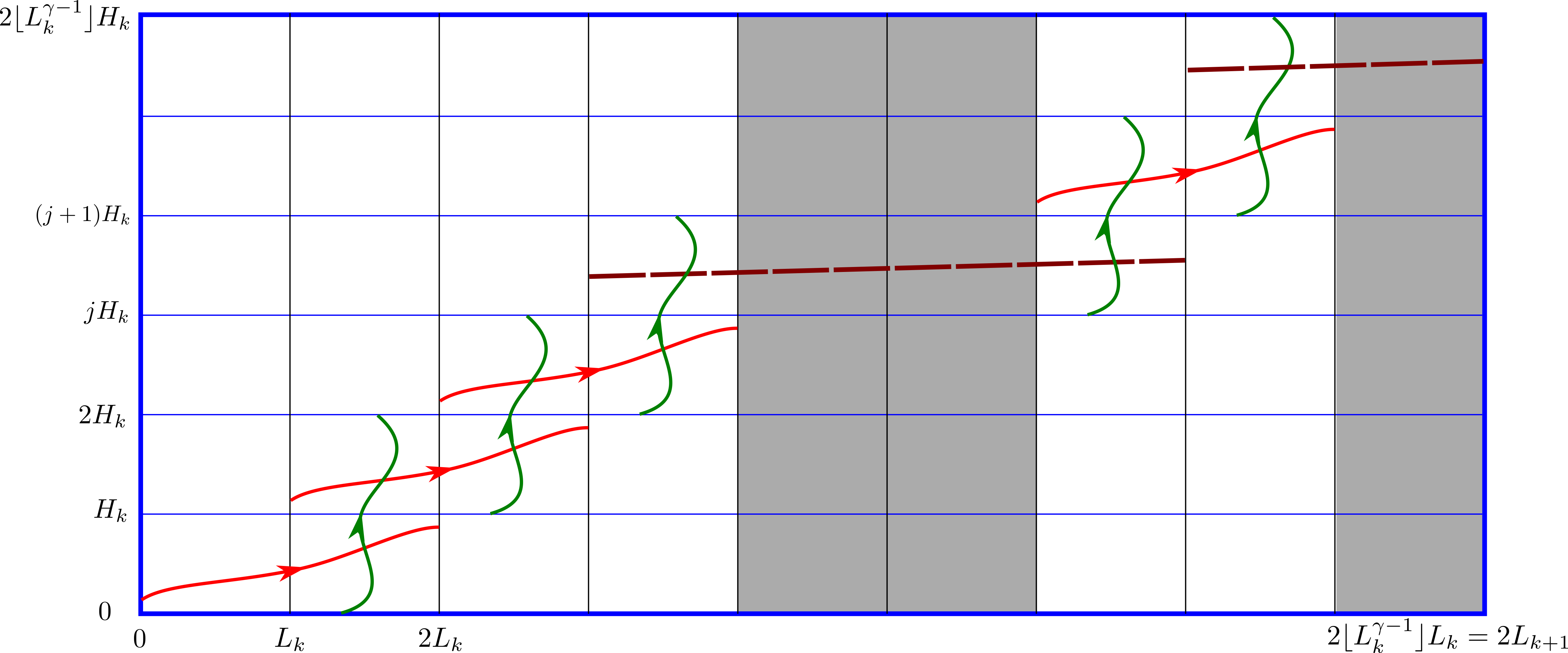}
    \caption{Each red line segment is an event of the form $\LRC(R^{\hor}_k(j,i))$ while the green ones signify the event $\DTC(R^{\ver}_k(j,i))$. By patching together horizontal and vertical crossings in an alternating manner, we obtain a left-right crossing of the wide rectangle. Bad areas (gray) are traversed differently (dotted red line), namely on the fastest path (see below). The total number of red and green crossings is at most  $2\cdot2\lfloor L_k^{\gamma-1}\rfloor$, each happening with probability at least $1-q_k(\lambda)$.}
    \label{fig:LRC}
\end{figure}

Let us specify next how bad areas are being traversed. Assume that $I_{k,i}$ is a bad scale-\(k\) block inside $I_{k+1,0}\cup I_{k+1,1}$ such that either $i=0$ or $I_{k,i-1}$ is a good block. Set $\underline{i} := (i-1)\vee 0$ and search for the first good block after $I_{k,i}$, which has index $\overline{i}:=\min\{\tilde{i}\geq i \colon  I_{k,\tilde{i}}\subset I_{k+1,0}\cup I_{k+1,1}\text{ is good}\}$ where we set $\overline{i}:=2L_k$ if the set is empty. Since $I_{k+1,0},I_{k+1,1}$ are good, there are no more than four bad blocks, hence $|\overline{i}+1-\underline{i}|\leq 6$. Let us construct a left-right crossing in a rectangle $[jH_k,(j+1)H_k]\times[\underline{i}L_k,\overline{i}L_k]$ for some $j$. We enumerate the points of  $\left(X\cap[\underline{i}L_k,\overline{i}L_k]\right)$ in increasing order as $(Z_x)_{x=0}^{m+1}$, denoting by \(m+1\) the total number of such points. By \eqref{eq:def:reduced-rectangle}, we do not consider the right-most column and is suffices to reach $Z_m$ for a successful left-right crossing. 
Now, we either have $(jH_k+x,Z_x)\in V$ or $(jH_k+x+1,Z_x)\in V$. We assume the former, otherwise replace all $jH_k+x$ terms accordingly.
A left-right crossing is guaranteed by opening all vertices 
\[
	(jH_k+x,Z_x),\quad 0\leq x \leq m
\]
and all edges
\[
	e_{x,x+1}:=(jH_k+x,Z_x)\to (jH_k+x+1,Z_{x+1}), \quad 0\leq x\leq m-1.
\]
As $H_k>6L_k$, this gives a crossing of the rectangle $[jH_k,(j+1)H_k]\times[\underline{i}L_k,\overline{i}L_k]$. Hence, the probability of traversing a bad area is no smaller than
\begin{equation}
   \begin{aligned}
       \P^X_\lambda\big(\LRC([jH_k,(j+1)H_k]    \times[\underline{i}L_k,\overline{i}L_k])\big)
       &\geq \prod_{x=0}^{m-1} \kappa_\lambda(\nu_{x,x+1})\prod_{x=0}^{m} \kappa_\lambda(\xi_{x})\geq \prod_{x=0}^{m} \e^{\nu_{x,x+1}}\e^{\xi_x}
       =\e^{Z_{m+1}-Z_0}
       \\ & \geq \e^{-6 L_k},
   \end{aligned} 
\end{equation}
using \eqref{eq:OPRE-stationarised-renewal-X} in the second and the definition of \(X\) in the last step, cf.~\eqref{eq:OPRE-requirements-connection-kernel}. Note that there are at most two bad areas that must be traversed. The remaining calculations follow those of~\cite{hilario2023phasetransition}. As crossing events are monotone and therefore positively correlated by the FKG inequality, we obtain
\begin{equation*}
    \begin{aligned}
        \P^X_\lambda(G_{k+1}) &\geq \P^X_\lambda(\text{all crossing events in Figure~\ref{fig:LRC} occur})\\
        & \geq \big(1-q_k(\lambda)\big)^{4\lfloor L_k^{\gamma-1}\rfloor}  \e^{-12 L_k}
        \geq \big(1-4\lfloor L_k^{\gamma-1}\rfloor q_k(\lambda)\big) \e^{-12 L_k}.
    \end{aligned}
\end{equation*}
By assumption, $q_k(\lambda)\leq \exp(-L_k^\beta)$ for large $k$. Also, by choosing $k$ larger than some $c_6=c_6(L_0,\gamma,\beta)$, we may assume $4\lfloor L_k^{\gamma-1}\rfloor \exp(-L_k^\beta) \leq 1-1/\e$ to infer
\begin{equation*} %\label{eq:estimate-G_k-horizontal-crossing}
    \P^X_\lambda(G_{k+1}) \geq %1/\e \cdot \e^{-12 L_k} \geq 
    \e^{-13 L_k}.
\end{equation*}
Recalling the estimates on \(L_k\) given in~\eqref{eq:estimate-L_k}, the claim follows hence as
\begin{equation*}
    \begin{aligned}
        \frac{\P^X_\lambda(\neg\LRC(R^\hor_{k+1}(0,0)))}{\exp(-L^\beta_{k+1})} 
        &
        	\leq \exp(L^\beta_{k+1})\P^X_\lambda(\neg G_{k+1})^{\lceil\exp(L^\mu_{k+1})\rceil}
       		\leq \exp(L^\beta_{k+1})\left( 1- \exp(-13L_k)\right)^{\exp(L^\mu_{k+1})}
       	\\ &
       		\leq \exp(L^\beta_{k+1}) \exp\left(-\exp(L^\mu_{k+1} - 13L_k))\right)
       		\leq \exp\big(  L_k^{\gamma\beta} - \exp(-13L_k + (1/2)^\mu L_k^{\gamma\mu})\big),
    \end{aligned}
\end{equation*}
which is smaller than $1$ for every $k$ larger than some $c_4=c_4(\gamma,L_0,\mu,\beta,c_6)$ since $\gamma\mu>1$.
\end{proof}

\begin{proof}[Proof of Lemma~\ref{lem:OPRE-induction-vertical-crossing}] 
    %We follow along the lines of the original proof and add some arguments from oriented percolation, e.g.~from \cite[Section 10]{durrett1984oriented}.\\
    Fix an environment $X$ in which $I_{k+1,0}$ is good. Write 
    \[
    	l:=\lfloor \tfrac12\lfloor L_k^{\gamma-1}\rfloor\rfloor\qquad\text{and}\qquad
    	T:=2H_{k+1}/H_k = 4 \lfloor L_{k+1}^{\gamma-1}\rfloor \lceil\exp(L_{k+1}^\mu)\rceil
   	\]
    and note that $l\geq2$ since $L_0^{\gamma-1}\geq5$. Hence, either all the $I_{k,i}$ with $i\in\{0,1,\dots,l-1\}$ are good or all those with $i\in\{l+1,\dots,\lfloor L_k^{\gamma-1}\rfloor\}$.  Without loss of generality, we consider the first case. The proof follows Peierls' argument in a renormalised lattice, see Figure~\ref{fig:vertical_crossing_renorm}. 
    
    \begin{figure}[ht]
        \centering
        \includegraphics[width=0.9\linewidth]{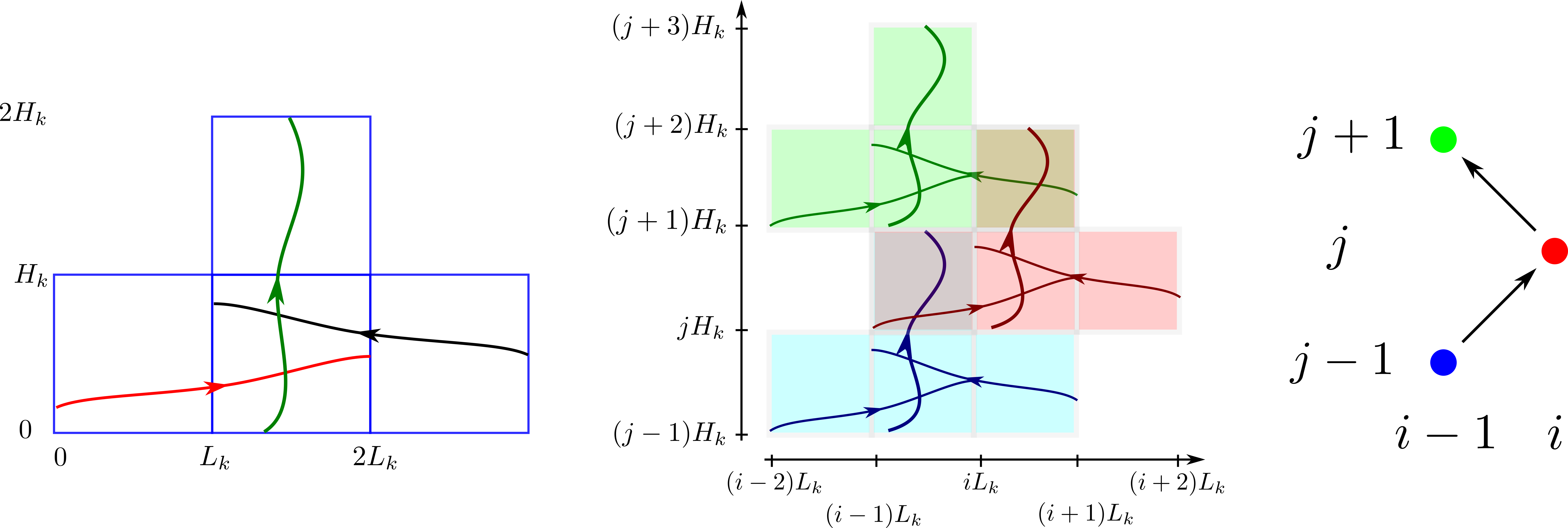}
        \caption{Left: A site is declared open if the corresponding left-right, right-left and bottom-top crossing exist. Middle: Open sites connect to each other as in the $\Ltwo$ graph. Right: We renormalise to the $\Ltwo$-lattice.}
        \label{fig:vertical_crossing_renorm}
    \end{figure}  
    
    For $j\in\{0,1,\dots,T\}$ and $i\in\{0,1,\dots,l\}$ with \(j+i\) even, we call a site $(j,i)$ \emph{open} if the three events $\LRC(R^\hor_k(j,i-1))$, $\RLC(R^\hor_k(j,i))$ and $\DTC(R^\ver_k(j,i))$ occur. (In the cases of $i-1<0$ and $i>l-1$, we ignore the corresponding event and consider it as having occurred.) Thus, a site is {\em closed} with probability no larger than
    \[
    	1-(1-q_k(\lambda))^3 \leq 3 q_k(\lambda) \leq 3\exp(-L_k^\beta),
    \]
    by assumption. This yields a $1$-dependent site percolation on $\Ltwo$. If  $\DTC(R^\ver_{k+1}(0,0))$ does not occur, then there is no bottom-top crossing in the rescaled lattice on $\Ltwo$ in the rectangle $[0,T]\times[0,l-1]=[0,4\lfloor L_{k+1}^{\gamma-1} \rfloor\lceil \exp(L_{k+1}^\mu)\rceil]\times [0,l-1]$ and there must be a blocking contour present that is of length exceeding $l-1$ and starts in $[0,T]\times \{0\}$, see Figure~\ref{fig:OPRE-contour-argument}.
    \begin{figure}[ht]
        \centering
        \includegraphics[width=0.5\linewidth]{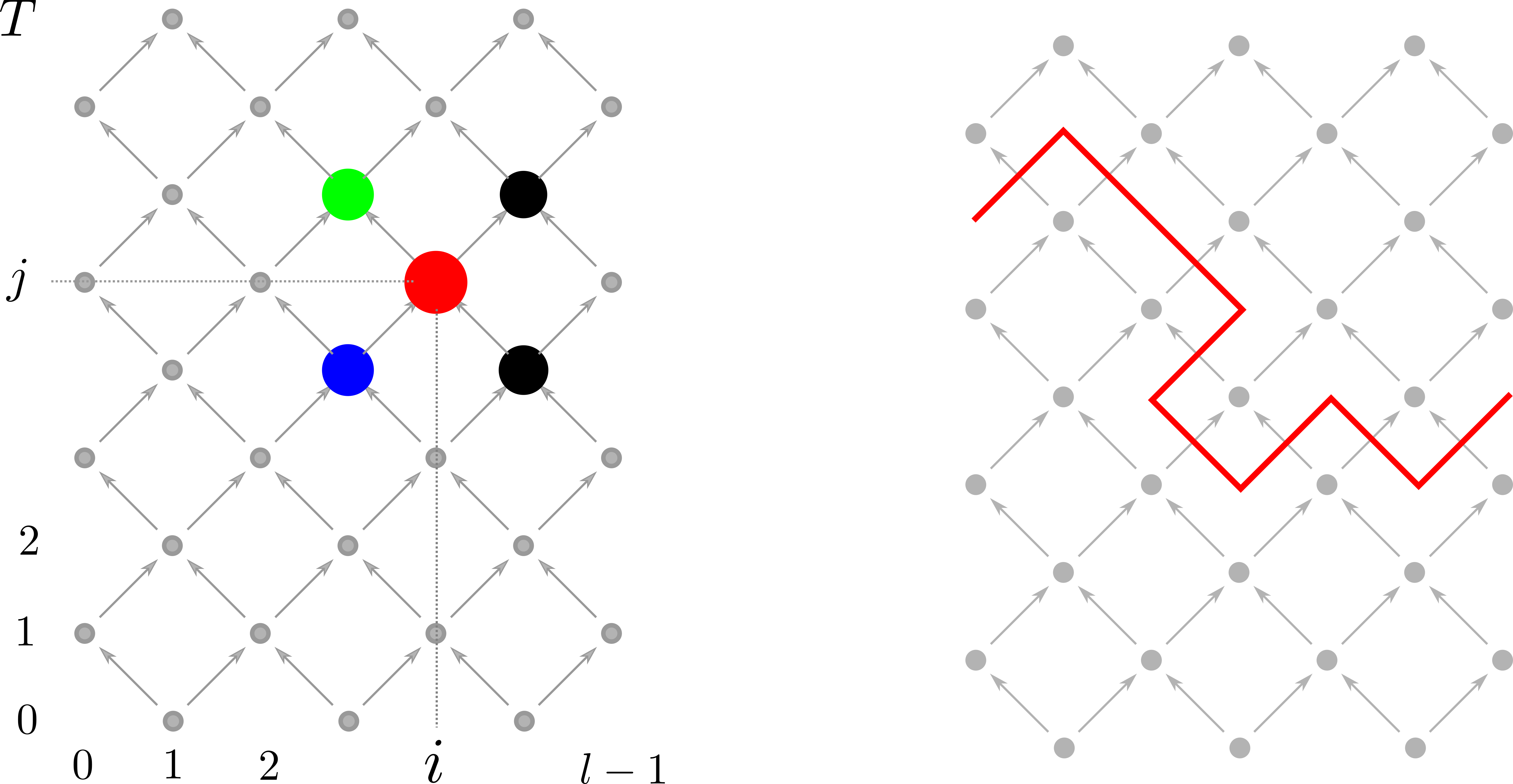}
        \caption{Renormalised subgraph. Left: The state of site $(j,i)$ depends only on its (at most) 4 neighbours. Right: If no bottom-top crossing exists, then there must be a blocking contour (red). Here, it has length $m=8$.}
        \label{fig:OPRE-contour-argument}
    \end{figure} 
    Note that a blocking contour of length $m$ can only exist, if there are no less than $m/4$ closed sites. Since the considered site percolation on \(\Ltwo\) is $1$-dependent, the state of a site $(t,x)$ depends on its four neighbours and itself. Hence, for each blocking contour, we can identify at least $\lceil m/20 \rceil$ many mutually independent sites that need to be closed. Note that there are no more than $4^m$-many contours in $\Ltwo$ starting at a fixed site. Therefore, we have
    \begin{equation*}%\label{eq:estimation-vertical-crossing}
        \begin{aligned}
            \P^X_\lambda(\neg \DTC & (R^{\ver}_{k+1}(0,0))) \leq
            \P^X_\lambda(\exists\, \text{blocking contour in renormalised lattice}) \\
            &\leq \sum_{m=l-1}^\infty \sum_{j=0}^{T} \P^X_\lambda\big(\exists\, \text{blocking contour in renormalised lattice of length } m \text{ starting in } (j,0)\big)\\
            &\leq \sum_{m=l-1}^\infty 4^m T \, \big(3\exp(-L_k^\beta)\big)^{m/20}
            \leq T  \sum_{m=l-1}^\infty \big(12\exp(-L_k^\beta/20)\big)^m\\
            &\leq 4\lfloor L_{k+1}^{\gamma-1} \rfloor\lceil \exp(L_{k+1}^\mu)\rceil \big(12\exp(-L_k^\beta/20)\big)^{l-1} c_7(\beta,L_0,\gamma)
            \leq c_7 \exp\big(L_{k+1}^\mu - c_8 L_k^{\beta + \gamma - 1}\big),
        \end{aligned}
    \end{equation*}

    for some $c_7=c_7(\gamma,L_0,\beta)$ sufficiently large, $c_8=c_8(\gamma,L_0,\beta)$ sufficiently small and all $k$ sufficiently large. The statement of Lemma~\ref{lem:OPRE-induction-vertical-crossing} then follows from
    \begin{equation*}
        \begin{aligned}
            \frac{\P^X_\lambda(\neg \DTC (R^{\ver}_{k+1}(0,0)))}{\exp(-L^\beta_{k+1})} \leq c_7 \exp(L_k^{\gamma\beta} + L_k^{\gamma\mu} - c_8 L_k^{\beta+\gamma-1}) \leq 1,
        \end{aligned}
    \end{equation*}
    for every $k\geq c_5=c_5(\gamma,L_0,\mu,\beta)$ sufficiently large, where we use \(\beta+\gamma-1>\max\{\gamma\beta,\gamma\mu\}\) in the last step.
\end{proof}

\subsection{Proof of main results} \label{subsec:proof-temporal-environment}
%We have now everything in place to prove our main theorem. 

\begin{proof}[Proof of Theorem~\ref{thm:OPRE-percolates}] 
%As always, the proof follows the same structure as the original \cite[Thm.~1.1]{hilario2023phasetransition}. 
We choose $\lambda$ large enough such that Lemma~\ref{lem:OPRE-main-induction} is satisfied. Consider the event 
\[
	A_k=\{\text{all scale-}k\text{ blocks inside } I_{k+1,0}\cup I_{k+1,1} \text{ are good}\}.
\]
Then, $\P(A_k)\geq 1-\lfloor L_k^{\gamma-1}\rfloor L_k^{-\alpha}\geq 1- L_k^{-\alpha/2}$ by Lemma~\ref{lem:OPRE-prob-bad-interval} and hence
\begin{equation}
    \begin{aligned}
        \sum_{k=0}^\infty \P(\neg A_k) 
        \leq \sum_{k=0}^\infty L_k^{-\alpha/2} <\infty,
    \end{aligned}
\end{equation}
as $L_k$ grows superexponentially. Therefore $\neg A_k$ only occurs finitely often almost surely. Fix a realisation of $X$ and let $K(X)\in\N$ be the smallest scale such that $A_k$ occurs for all $k\geq K(X)$. 

\begin{figure}
    \centering
    \includegraphics[width=0.8\linewidth]{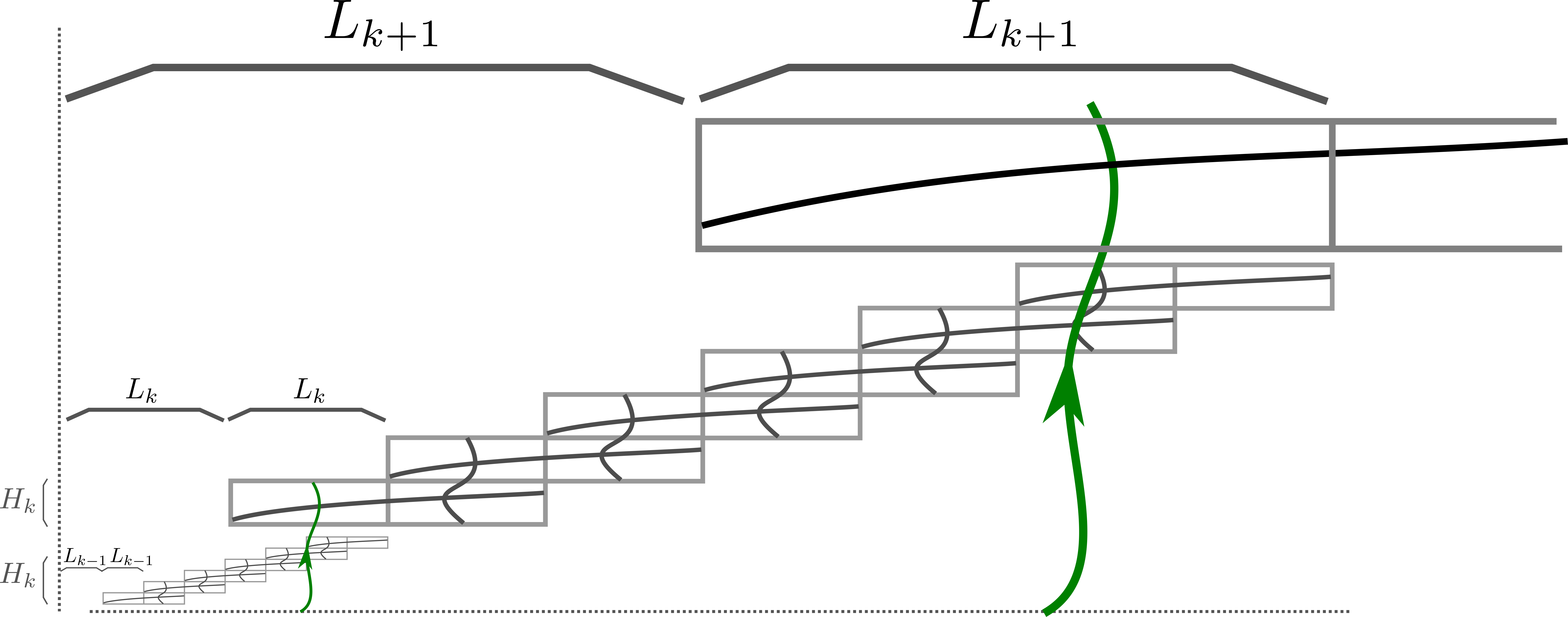}
    \caption{Construction of an infinite directed path based on the events $\mathcal{V}_k$ (green arrows) and $\mathcal{H}_k$ (black ragged staircase). Due to planarity, the paths intersect to form an infinite path.}
    \label{fig:construction_infinite_cluster}
\end{figure}

We build an infinite cluster by patching together left-right and bottom-top crossings of rectangles (as in~\eqref{eq:OPRE-def-hor-rectangle} and~\eqref{eq:def-ver-rectangle}) in an alternating manner, see Figure~\ref{fig:construction_infinite_cluster}. More precisely, we consider the events
$\mathcal{V}_k := \DTC(R_k^\ver(0,1))$
as well as
\[
\mathcal{H}_k:=\Big(\bigcup_{i=1}^{2\lfloor L_k^{\gamma-1}\rfloor-1}\LRC(R^\ver_k(i,i))\cup\DTC(R^\hor_k(i,i+1))\Big).
\]
Clearly, if $\mathcal{V}_k, \mathcal{H}_k$ only fail to occur for finitely many \(k\), an infinite connected path exists, see Figure~\ref{fig:construction_infinite_cluster}. However, this is almost surely the case as
\begin{equation*}
    \begin{aligned}
        \sum_{k=K(X)}^\infty \big(\P^X_\lambda(\neg \mathcal{V}_k)+\P^X_\lambda(\neg\mathcal{H}_k)\big)
        \leq \sum_{k=K(X)}^\infty \big(q_k + 4 \lfloor L_k^{\gamma-1}\rfloor q_k \big)
        \leq 5 \sum_{k=K(X)}^\infty L_k^{\gamma-1}\exp(-L_k^\beta) < \infty.
    \end{aligned}
\end{equation*}
This concludes the proof.
\end{proof}

%\subsection{Random temporal environment}\label{subsec:proof-temporal-environment}

It remains to prove Proposition~\ref{prop:OPRE-temporal-extinction} regarding the oriented percolation model with temporal stretches.

\begin{proof}[Proof of Proposition~\ref{prop:OPRE-temporal-extinction}]
    The number of vertices reached from the origin grows at most polynomially in time, so the idea is to find some $\nu_{t,t+1}$ sufficiently large such that all edges at time $t$ are closed. Let $p\in(0,1)$. At time $t$, there are at most $2t$ outgoing edges from vertices of the form $(t,x)$ with $x=0,\dots,t$. Given a realisation $\bar\nu:=(\nu_{t,t+1})_{t\in\Zz}$ of the temporal stretches, we estimate
    \begin{equation*}
        \begin{aligned}
            \P(o\leadsto\infty| \bar\nu) &\leq \prod_{t\in\Zz} \big(1-\P\big((t,x)\to(t+1,x+\sigma)\text{ is closed }\forall x\leq t+1,\,\sigma\in\{-1,1\}| \bar\nu\big)\big) \\
            &\leq \inf_{t\in\Zz} \big(1-\P\big((t,x)\to(t+1,x+\sigma)\text{ is closed }\forall x\leq t+1,\,\sigma\in\{-1,1\} |\bar\nu\big)\big)\\
            &= \inf_{t\in\Z_{\ge 1}} \left(1-(1-p^{\nu_{t-1,t}})^{2t}\right).
        \end{aligned}
    \end{equation*}
    It suffices to show $1 = \sup_{t\in\Z_{\ge 1}} (1-p^{\nu_{t-1,t}})^{2t}$ for almost-every realisation of the $\bar\nu$, or equivalently
    \begin{equation}\label{eq:proof:OPRE-temporal-extinction}
        \inf_{t\in\Z_{\ge 1}} \nu_{t-1,t}\log(p)+\log t = -\infty.
    \end{equation}
    Take $c_1=c_1(p)>0$ such that $c_1\log(p)<-1$. Then, \eqref{eq:proof:OPRE-temporal-extinction} is satisfied if $\nu_{t-1,t}\geq c_1\log t$ for infinitely many $t$. We will see that this is the case, thereby finishing the proof. Given some $c>1$,  \eqref{eq:heavy-tail-stretches} yields 
    $$\P(\nu_{t-1,t} > c_1 \log t) \geq c^{-c_1 \log t} = t ^ {- c_1 \log c}$$
    for infinitely many $t\in\Zz$. We enumerate those as $(t_k)_{k\in\Zz}\subset\Zz$. (This subsequence depends on $c_1=c_1(p)$.) Choosing $c>1$ such that $\delta := c_1\log c \leq 1$, we see that
    \begin{equation}
        \sum_{t=1}^\infty \P(\nu_{t-1,t} \geq c_1 \log t) \geq \sum_{k=1}^\infty \sum_{t=t_{k-1}+1}^{t_k} \P(\nu_{t-1,t} \geq c_1 \log t_k) \geq \sum_{k=1}^\infty (t_k - t_{k-1}) t_k^{-\delta} = \infty,
    \end{equation}
    since $\delta\leq1$. As the events \{$\nu_{t-1,t} \geq c_1 \log t$\} are independent from each other, the Borel--Cantelli lemma tells us that they happen infinitely often. Thus, the claim follows.
\end{proof}

\paragraph{Acknowledgement.} This research was supported by the Leibniz Association within the Leibniz Junior Research Group on {\em Probabilistic Methods for Dynamic Communication Networks} as part of the Leibniz Competition (grant no.\ J105/2020) and the Berlin Cluster of Excellence {\em MATH+} through the project {\em EF45-3} on {\em Data Transmission in Dynamical Random Networks}.

%%%%%%%%%%%%%%%%%%%%%%%
%%%%%Bibliography%%%%%%
%%%%%%%%%%%%%%%%%%%%%%%
\section*{References}
\renewcommand*{\bibfont}{\footnotesize}
\printbibliography[heading = none]
\end{document}